\documentclass[12pt,a4paper,reqno]{amsart}

\usepackage{amssymb}
\usepackage{amsmath}
\usepackage{latexsym}
\usepackage{exscale}
\usepackage[latin1]{inputenc}
\usepackage{epsfig}
\usepackage{verbatim}
\usepackage{graphicx}
\usepackage{subfigure}

\newcommand{\norm}[1]{\left\Vert#1\right\Vert}
\newcommand{\abs}[1]{\left| #1 \right|}
\newcommand{\ip}[2]{\langle #1 , #2 \rangle}

\numberwithin{equation}{section} 

\newtheorem{thm}{Theorem}[subsection] 

\newtheorem{defi}[thm]{Definition}
\newtheorem{rem}[thm]{Remark}

\newtheorem{lem}[thm]{Lemma}

\begin{document}

\allowdisplaybreaks

\title[Triebel-Lizorkin Spaces and Shearlets] {\large
Triebel-Lizorkin Spaces and Shearlets on the Cone in $\mathbb{R}^2$}

\author{Daniel Vera}

\address{Daniel Vera
\\
Departamento de Matem\'aticas
\\
Universidad Aut\'onoma de Madrid
\\
28049 Madrid, Spain}

\email{daniel.vera@uam.es}

\date{\today}


\begin{abstract} The shearlets are a special case of the wavelets with composite
dilation that, among other things, have a basis-like structure and
multi resolution analysis properties. These relatively new
representation systems have encountered wide range of applications,
generally surpassing the performance of their ancestors due to their
directional sensitivity. However, little is known about their
relation with spaces other than $L^2$. Here, we find a
characterization of a kind of anisotropic inhomogeneous
Triebel-Lizorkin spaces (to be defined) with the so called
``shearlets on the cone" coefficients. We first prove the
boundedness of the analysis and synthesis operators with the
``traditional" shearlets coefficients. Then, with the development of
the smooth Parseval frames of shearlets of Guo and Labate we are
able to prove a reproducing identity, which was previously possible
only for the $L^2$ case. We also find some embeddings of the
(classical) dyadic spaces into these highly anisotropic spaces, and
viceversa, for certain ranges of parameters. In order to keep a
concise document we develop our results in the ``weightless" case
($w=1$) and give hints on how to develop the weighted case.


\end{abstract}

\maketitle

\section{Introduction.}\label{S:intro}
The traditional (separable) multidimensional wavelets are built from
tensor-like products of $1$-dimensional wavelets. Hence, wavelets
cannot ``sense" the geometry of lower dimension discontinuities. In
$\mathbb{R}^d$ the number of wavelets are $2^d-1$ for each scale. In
applications it may be desirable to be able to detect more
orientations having still a basis-like representation. In recent
years there have been attempts to achieve this sensitivity to more
orientations. Some of them include the directional wavelets or
filterbanks \cite{AMV}, \cite{BS}, the \emph{curvelets}
\cite{CaDo00} and the \emph{contourlets} \cite{DoVe}, to name just a
few. Regarding the contourlets, since they are built in a
discrete-time setting and from a finite set of parameters, they lack
of flexibility and, for some applications, assume there exist smooth
spatially compactly supported functions approximating a frequency
partition as that in Subsection \ref{sS:Discrt_Shrlts_Cone}. On the
other hand, the curvelets are built on polar coordinates so their
implementation is rather difficult.

In \cite{GLLWW}, Guo, Lim, Labate, Weiss and Wilson, introduced the
wavelets with composite dilation. This type of representation takes
full advantage of the theory of affine systems on $\mathbb{R}^n$ and
therefore provides a natural transition from the continuous
representation to the discrete (basis-like) setting (as in the case
of wavelets). A special case of the composite dilation wavelets is
that of the shearlets system which provides Parseval frames for
$L^2(\mathbb{R}^2)$ or subspaces of it (depending on the discrete
sampling of parameters, see Subsections \ref{sS:Discrt_Shrlts} and
\ref{sS:Discrt_Shrlts_Cone}). A large amount of applications of the
shearlet transform not only to the image processing can be consulted
in \texttt{http://www.shearlet.org}.

We characterize a new kind of \emph{highly anisotropic
inhomogeneous} Triebel-Lizorkin spaces using the ``shearlets on the
cone" coefficients. The line of argumentation follows the
$\varphi$-transform in \cite{FJ90}.

Following the classical definition of Triebel-Lizorkin spaces by
Triebel \cite{Tr83}, \cite{Tr92}, Frazier-Jawerth \cite{FJ90},
Frazier-Jawerth-Weiss \cite{FJW} and their weighted counterparts in
the work of Bui \cite{Bu82}, \cite{Bu84}; more recently, Bownik and
Ho in \cite{BH05} define the \emph{weighted anisotropic
inhomogeneous Triebel-Lizorkin spaces}
$\mathbf{F}^{\alpha,q}_p(A,w)$ as the collection of all
$f\in\mathcal{S}'$ such that
$$\norm{f}_{\mathbf{F}^{\alpha,q}_p(A,w)}
    =\norm{f\ast\Phi}_{L^p(w)}+
    \norm{\left(\sum_{j=1}^\infty(\abs{\text{det}A}^{j\alpha}\abs{f\ast\varphi_j})^q\right)^{1/q}}_{L^p(w)}<\infty,$$
where $\Phi, \varphi\in\mathcal{S}$ with the properties that
$\text{supp}(\hat{\Phi})$ is compact and
$\text{supp}(\hat{\varphi})$ is compact and bounded away from $0$
and where $A\in GL_d(\mathbb{R})$ with all of its eigenvalues $>1$.
Then, they show that there exists another pair
$\Psi,\psi\in\mathcal{S}$ with the same properties such that
$$\overline{\hat{\Psi}(\xi)}\hat{\Phi}(\xi) + \sum_{j=1}^\infty \overline{\hat{\varphi}(\xi A^{-j})} \hat{\psi}(\xi A^{-j})
    =1,\;\;\; \text{for all }\xi\in\mathbb{\hat{R}}^d,$$
which yields the representation formula
$$f=\sum_{\abs{Q}=1}\ip{f}{\Phi_Q}\Psi_Q
    +\sum_{\abs{Q}<1}\ip{f}{\varphi_Q}\psi_Q,$$
for any $f\in\mathcal{S}'$ with convergence in $\mathcal{S}'$ and
where $Q$ runs through the ``cubes" $Q_{j,k}=A^{-j}((0,1]^d+k)$
(notice there is no shear operation as in (\ref{e:def_Qjlk})). This
obviates (ignores) the fact that for dimension $d$ one needs $2^d-1$
wavelets to cover $\mathbb{\hat{R}}^d$ (or $\mathbb{R}^d$). Since
the number of wavelets remains the same across scales, one can
ignore the sum over the set of non-zero vertices of the cube which
changes the norm only by a constant. This is not the case for the
shearlets since the cardinality of the shear parameter $\ell$ grows
with $j$ as $\ell=-2^j,...,2^j$. An observation from \cite{BH05} is
that ``in the standard dyadic case $A=2I$, where $I$ is the identity
matrix, and then the factor $\abs{\text{det
}A}^{j\alpha}=2^{nj\alpha}$ in the above definition, instead of the
usual $2^{j\alpha}$. Then, there is a re-scaling of the smoothness
parameter $\alpha$, which in the traditional case is thought of as
the number of derivatives." The same happens in the setting of the
``shearlets on the cone", as we will see.

Up to now, at least to our knowledge, there has been only one
attempt to relate the shearlets with spaces other than $L^2$ as done
by Dahlke, Kutyniok, Steidl and Teschke in \cite{DKST}. They
establish new families of smoothness spaces by means of the coorbit
space theory.

The method used here can be applied to the case of spatially compact
support (separable) shearlets since they are frames with irregular
sampling (see \cite{KiKuLi}). This method cannot be applied to the
case of the discrete shearlets (see Subsection
\ref{sS:Discrt_Shrlts}) since the shear parameter $\ell\in
\mathbb{Z}$ causes that Lemma \ref{l:Nested_rectngls_ellipsoids}
fails, avoiding a characterization of the homogeneous case.
Nevertheless, this method can be applied to higher dimensions and
different anisotropic and shear matrices, as long as a kind of Lemma
\ref{l:Nested_rectngls_ellipsoids} holds.


The outline of the paper is as follows. We review the basic facts of
the different shearlet transforms in Section \ref{S:Shrlts} and give
the pertinent results with the corresponding references. In Section
\ref{S:Notation} we set notation and give two basic lemmata
regarding almost orthogonality in the ``shearlets on the cone"
setting. In Sections \ref{S:T-L-Shearlets} and
\ref{S:Ident_SmthPrsvlFrms} we mainly follow \cite{FJ90},
\cite{FJ85}, \cite{FJW} and \cite{BH05} to 1) prove the
characterization in terms of the ``shearlets on the cone"
coefficients and 2) prove the identity on $\mathcal{S}'$. In Section
\ref{S:Embeddings} we prove some embbeddings of (classical) dyadic
inhomogeneous Triebel-Lizorkin spaces into the highly anisotropic
inhomogeneous Triebel-Lizorkin spaces, and viceversa, for a certain
range of the smoothness parameter. In Section \ref{S:Wghts} we
explain how to extend this work to the weighted case. Proofs for
Sections \ref{S:Notation} and \ref{S:T-L-Shearlets} are given in
Section \ref{S:Proofs_lemmata}.

\vskip1cm
\section{Shearlets}\label{S:Shrlts}
The shearlets are a generalization of the wavelets which better
capture the geometrical properties of functions. They are also a
special case of the so-called wavelets with composite dilation (see
\cite{GLLWW}). We give a basic introduction to the construction of
different type of shearlets: continuous, discrete and discrete on
the cone, in the next three subsections. A point $x\in\mathbb{R}^d$
is a column vector $x=(x_1,\ldots,x_d)^t$ and a point $\omega$ in
the dual $\hat{\mathbb{R}}^d$ is a row vector
$\omega=(\omega_1,\ldots,\omega_d)$.

\vskip0.5cm
\subsection{Continuous shearlets}\label{sS:Cont_Shrlts}
A \textbf{continuous affine system} in $L^2(\mathbb{R}^d)$ is a
collection of functions of the form
$$\{T_tD_M\psi: t\in\mathbb{R}^d, M\in G\},$$
where $\psi\in L^2(\mathbb{R}^d)$, $T_t$ is the translation operator
$T_tf(x)=f(x-t)$, $D_M$ is the dilation operator
$D_Mf(x)=\abs{\text{det} M}^{-1/2}f(M^{-1}x)$ (normalized in
$L^2(\mathbb{R}^d)$), and $G$ is a subset of $GL_d(\mathbb{R})$. In
the case $d=2$, $G$ is the $2$-parameter dilation group (see
\cite{GLLWW} for an even more general definition)
$$G=\{M_{as}= \left(%
\begin{array}{cc}
  a & \sqrt{a}s \\
  0 & \sqrt{a} \\
\end{array}%
\right) :(a,s)\in\mathbb{R}_+\times \mathbb{R}\}.$$ The matrix
$M_{as}$ is the product $S_sA_a$ where $S_s= \left(%
\begin{array}{cc}
  1 & s \\
  0 & 1 \\
\end{array}%
\right)$ is the area preserving \textbf{shear transformation}
and $A_a=\left(%
\begin{array}{cc}
  a & 0 \\
  0 & \sqrt{a} \\
\end{array}%
\right)$ is the \textbf{anisotropic dilation}. Assume, in addition,
that $\psi$ is given by
\begin{equation}\label{e:Shrlts_Composition}
\hat{\psi}(\xi)=\hat{\psi}_1(\xi_1)\hat{\psi}_2(\frac{\xi_2}{\xi_1}),
\end{equation}
for any $\xi=(\xi_1,\xi_2)\in\hat{\mathbb{R}}^2$, $\xi_1\not=0$, and
where $\psi_1$ satisfies (Calderón's admissibility condition)
$$\int_0^\infty \abs{\hat{\psi}_1(a\omega)}^2\frac{da}{a}=1, \;\;\; \text{for a.e. } \omega\in\mathbb{R},$$
and $\norm{\psi_2}_{L^2(\mathbb{R})}=1$. Then, the affine system
$$\{\psi_{ast}(x)=a^{-3/4}\psi(M_{as}^{-1}(x-t)):a\in\mathbb{R}_+,s\in\mathbb{R}, t\in\mathbb{R}^2\},$$
is a reproducing system for $L^2(\mathbb{R}^2)$, that is,
$$\norm{f}_{L^2(\mathbb{R}^2)}^2=\int_{\mathbb{R}^2}\int_{\mathbb{R}}\int_0^\infty \abs{\ip{f}{\psi_{ast}}}^2 \frac{da}{a^3}ds dt,$$
for all $f\in L^2(\mathbb{R}^2)$ (see \cite{WW01}). 

\vskip0.5cm
\subsection{Discrete shearlets}\label{sS:Discrt_Shrlts}
Since $L^2(\mathbb{R}^2)$ is a separable Hilbert space, it happens
that, by an appropriate ``sampling" of the parameters of the
continuous shearlets, there exists a construction of a basis-like
system for $L^2(\mathbb{R}^2)$ with a ``discrete" shearlet system.

A countable family $\{e_j:j\in \mathcal{J}\}$ of elements in a
separable Hilbert space $\mathbb{H}$ is called a \textbf{frame} if
there exist constants $0<A\leq B<\infty$, such that
$A\norm{f}_{\mathbb{H}}^2\leq\sum_{j\in
\mathcal{J}}\abs{\ip{f}{e_j}}^2\leq B\norm{f}_{\mathbb{H}}^2$, for
all $f\in \mathbb{H}$. A frame is called \textbf{tight} if $A=B$,
and is called a \textbf{Parseval frame} if $A=B=1$. Thus, if
$\{e_j:j\in \mathcal{J}\}$ is a Parseval frame for $\mathbb{H}$,
then
$\norm{f}_\mathbb{H}^2=\sum_{j\in\mathcal{J}}\abs{\ip{f}{e_j}}^2$,
for all $f\in\mathbb{H}$, which is equivalent to the reproducing
formula $f=\sum_{j\in\mathcal{J}} \ip{f}{e_j}e_j$, with
convergence in $\mathbb{H}$.

%

With a special sampling of the parameters one can construct a
Parseval frame of discrete shearlets for $L^2(\mathbb{R}^2)$ (see
\cite{GLLWW}).

\vskip0.5cm
\subsection{Discrete shearlets on the cone}\label{sS:Discrt_Shrlts_Cone}
In spite of the Parseval frame property for $L^2(\mathbb{R}^2)$ of
the discrete shearlets system, there are no ``equivalent
information" among the ``mostly horizontal" and ``mostly vertical"
shearlets  (important in applications) since the tiling of
$\hat{\mathbb{R}}^2$ is not ``homogeneous" in these directions. The
covering of $\hat{\mathbb{R}}^2$ by the discrete shearlets is done
firstly by vertical bands or strips related to mostly horizontal
dilations (indexed by $j$). Then, each band is covered by infinitely
countable shear (area preserving) transformations (indexed by
$\ell$). With a little modification on the discrete shearlet system
above one can obtain a Parseval frame for functions in
$L^2(\mathbb{R}^2)$ whose Fourier transform is supported in the
horizontal cone
\begin{equation}\label{e:Cone_Domain}
\mathcal{D}^h=\{(\xi_1,\xi_2)\in\hat{\mathbb{R}}^2:
    \abs{\xi_1}\geq\frac{1}{8}, \abs{\frac{\xi_2}{\xi_1}}\leq 1\}.
\end{equation}
Let now $\hat{\psi}_1,\hat{\psi}_2\in C^\infty(\mathbb{R})$ with
$\text{supp }\hat{\psi}_1\subset
[-\frac{1}{2},-\frac{1}{16}]\cup[\frac{1}{16},\frac{1}{2}]$ and
$\text{supp }\hat{\psi}_2\subset [-1,1]$ such that
\begin{equation}\label{e:Discrt_Shrlt_Cond_Cone_1}
\sum_{j\geq 0}\abs{\hat{\psi}_1(2^{-2j}\omega)}^2 =1, \;\;\;
\text{for }\abs{\omega}\geq \frac{1}{8}
\end{equation}
and
\begin{equation}\label{e:Discrt_Shrlt_Cond_Cone_2}
\abs{\hat{\psi}_2(\omega-1)}^2+\abs{\hat{\psi}_2(\omega)}^2+\abs{\hat{\psi}_2(\omega+1)}^2=1,
\;\;\; \text{for } \abs{\omega}\leq 1.
\end{equation}
It follows from (\ref{e:Discrt_Shrlt_Cond_Cone_2}) that, for
$j\geq 0$,
\begin{equation}\label{e:Discrt_Shrlt_Cond_Cone_3}
\sum_{\ell=-2^j}^{2^j} \abs{\hat{\psi}_2(2^j\omega-\ell)}^2=1,
\;\;\; \text{for }\abs{\omega}\leq 1.
\end{equation}
Let
$$A_h= \left(%
\begin{array}{cc}
  4 & 0 \\
  0 & 2 \\
\end{array}%
\right), \;\;\; B_h=\left(%
\begin{array}{cc}
  1 & 1 \\
  0 & 1 \\
\end{array}%
\right)$$ and
$\hat{\psi}^h(\xi)=\hat{\psi}_1(\xi_1)\hat{\psi}_2(\frac{\xi_2}{\xi_1})$.
From (\ref{e:Discrt_Shrlt_Cond_Cone_1}) and
(\ref{e:Discrt_Shrlt_Cond_Cone_3}) it follows that
\begin{eqnarray}\label{e:PrsvlFrm_Prop_Shrlts}
\nonumber
  \sum_{j\geq 0}\sum_{\ell=-2^j}^{2^j} \abs{\hat{\psi}^h(\xi A^{-j}_hB^{-\ell}_h)}^2
    &=& \sum_{j\geq 0}\sum_{\ell=-2^j}^{2^j} \abs{\hat{\psi}_1(2^{-2j}\xi_1)}^2
        \abs{\hat{\psi}_2(2^j\frac{\xi_2}{\xi_1}-\ell)}^2 \\
    &=& \sum_{j\geq 0} \abs{\hat{\psi}_1(2^{-2j}\xi_1)}^2
        \sum_{\ell=-2^j}^{2^j}
        \abs{\hat{\psi}_2(2^j\frac{\xi_2}{\xi_1}-\ell)}^2=1,
\end{eqnarray}
for $\xi=(\xi_1,\xi_2)\in\mathcal{D}^h$ and which we will call the
\textbf{Parseval frame condition} (for the horizontal cone). Since
$\text{supp }\hat{\psi}^h\subset [-\frac{1}{2},\frac{1}{2}]^2$,
(\ref{e:PrsvlFrm_Prop_Shrlts}) implies that the shearlet system
\begin{equation}\label{e:Shrlt_Sys_Cone}
\{\psi_{j,\ell,k}^h(x)= 2^{3j/2}\psi^h(B^\ell_h A^j_h x-k): j\geq
0, -2^j\leq\ell\leq2^j, k\in\mathbb{Z}^2\},
\end{equation}
is a Parseval frame for $L^2((\mathcal{D}^h)^\vee)=\{f\in
L^2(\mathbb{R}^2): \text{supp }\hat{f}\subset\mathcal{D}^h\}$ (see
\cite{GLLWW}, Subsection 5.2.1). This means that
$$\sum_{j\geq 0}\sum_{\ell=-2^j}^{2^j}\sum_{k\in\mathbb{Z}^2}
    \abs{\ip{f}{\psi_{j,\ell,k}^h}}^2 = \norm{f}^2_{L^2(\mathbb{R}^2)},$$
for all $f\in L^2(\mathbb{R}^2)$ such that $\text{supp
}\hat{f}\subset \mathcal{D}^h$. There are several examples of
functions $\psi_1,\psi_2$ satisfying the properties described above
(see \cite{GL07}). Since $\hat{\psi}^h\in
C^\infty_c(\hat{\mathbb{R}}^2)$,  there exist $C_N$ such that
$\abs{\psi^h(x)}\leq C_N (1+\abs{x})^{-N}$ for all $N\in\mathbb{N}$.
The geometric properties of the horizontal shearlets system are more
evident by observing that
$$\text{supp }(\psi_{j,\ell,k})^\wedge\subset \{\xi\in\hat{\mathbb{R}}^2
: \xi_1\in [-2^{2j-1}, -2^{2j-4}]\cup[2^{2j-4},2^{2j-1}],
\abs{\frac{\xi_2}{\xi_1}-\ell 2^{-j}}\leq 2^{-j}\}.$$

One can also construct a Parseval frame for the vertical cone
\begin{equation*}
\mathcal{D}^v=\{(\xi_1,\xi_2)\in\hat{\mathbb{R}}^2:
    \abs{\xi_2}\geq\frac{1}{8}, \abs{\frac{\xi_1}{\xi_2}}\leq 1\},
\end{equation*}
by defining
$\hat{\psi}^v(\xi)=\hat{\psi}_1(\xi_2)\hat{\psi}_2(\frac{\xi_1}{\xi_2})$
and with anisotropic and shear matrices
$$A_v= \left(%
\begin{array}{cc}
  2 & 0 \\
  0 & 4 \\
\end{array}%
\right), \;\;\; B_v=\left(%
\begin{array}{cc}
  1 & 0 \\
  1 & 1 \\
\end{array}%
\right).$$

Let $\hat{\varphi}\in C^\infty_c(\mathbb{R}^2)$, with $\text{supp
}\hat{\varphi}\subset [-\frac{1}{4},\frac{1}{4}]^2$ and
$\abs{\hat{\varphi}}=1$ for $\xi\in
[-\frac{1}{8},\frac{1}{8}]^2=\mathcal{R}$, be such that
\begin{eqnarray}\label{e:Shrlt_Repr_Sys}
\nonumber
  P(\xi)
    &=& \abs{\hat{\varphi}(\xi)}^2\chi_\mathcal{R}(\xi)
    + \sum_{j\geq 0}\sum_{\ell=-2^j}^{2^j}\abs{\hat{\psi}^h(\xi A^{-j}_hB^{-\ell}_h)}^2\chi_{\mathcal{D}^h}(\xi) \\
    &&\;\;\;\; +\sum_{j\geq 0}\sum_{\ell=-2^j}^{2^j}\abs{\hat{\psi}^v(\xi
        A^{-j}_vB^{-\ell}_v)}^2\chi_{\mathcal{D}^v}(\xi)=1,
        \;\;\;\text{for all } \xi\in\hat{\mathbb{R}}^2.
\end{eqnarray}

\vskip1cm
\section{Notation and almost orthogonality}\label{S:Notation}
Since all results in the horizontal cone $\mathcal{D}^h$ can be
stated for the vertical one $\mathcal{D}^v$, with the obvious
modifications as explained in Subsection
\ref{sS:Discrt_Shrlts_Cone}, we drop the superindex $h$ and develop
only for the horizontal cone and refer only to ``the cone".

We will develop our results with
$$A= \left(%
\begin{array}{cc}
  4 & 0 \\
  0 & 2 \\
\end{array}%
\right) \;\;\; \text{ and }\;\;\; B=\left(%
\begin{array}{cc}
  1 & 1 \\
  0 & 1 \\
\end{array}%
\right),$$ an anisotropic dilation and a shear matrix, respectively.
We consider $\psi$ defined by $\hat{\psi}(\xi)=
\hat{\psi}_1(\xi_1)\hat{\psi}_2(\xi_2/\xi_1)$ with $\psi_1$ and
$\psi_2$ satisfying (\ref{e:Discrt_Shrlt_Cond_Cone_1}),
(\ref{e:Discrt_Shrlt_Cond_Cone_2}) and
(\ref{e:Discrt_Shrlt_Cond_Cone_3}). In order to follow \cite{FJ90}
we will require $\varphi, \psi\in\mathcal{S}$ with the same
conditions on $\text{supp }(\hat{\varphi})$, $\text{supp
}(\hat{\psi}_1)$, $\text{supp }(\hat{\psi}_2)$ and equations
(\ref{e:Discrt_Shrlt_Cond_Cone_1}),
(\ref{e:Discrt_Shrlt_Cond_Cone_2}) and
(\ref{e:Discrt_Shrlt_Cond_Cone_3}) in order to preserve the
geometrical properties of anisotropic and shear operations.
Following the notation for the usual isotropic dilation,
$\varphi_t(x):=\frac{1}{t}\varphi(\frac{x}{t})$, we denote for a
matrix $M\in GL_2(\mathbb{R})$ the anisotropic dilation
$\varphi_M(x)=\abs{\text{det} M}^{-1}\varphi(M^{-1}x)$ (do not
confuse with the dilation operator normalized in $L^2$ as in
Subsection \ref{sS:Cont_Shrlts}). We also denote
$\tilde{\varphi}(x)=\overline{\varphi(-x)}$. For $Q_0=[0,1)^2$,
write
\begin{equation}\label{e:def_Qjlk}
Q_{j,\ell,k}=A^{-j}B^{-\ell}(Q_0+k),
\end{equation}
with $j\geq 0$, $\ell=-2^j,\ldots, 2^j$ and $k\in\mathbb{Z}^2$.
Therefore, $\int
\chi_{Q_{j,\ell,k}}=\abs{Q_{j,\ell,k}}=\abs{Q_{j,\ell}}=\abs{Q_j}=2^{-3j}=\abs{\text{det
} A}^{-j}$. We also write
$\tilde{\chi}_Q(x)=\abs{Q}^{-1/2}\chi_Q(x)$. Let
$\mathcal{Q}_{AB}:=\{Q_{j,\ell,k}: j\geq 0,\ell=-2^j,\ldots, 2^j,
k\in\mathbb{Z}^2\}$ and $\mathcal{Q}^{j,\ell}:=
\{Q_{j,\ell,k}:k\in\mathbb{Z}^2\}$, then $\mathcal{Q}^{j,\ell}$ is a
partition of $\mathbb{R}^2$. To shorten notation and clear
exposition, we will identify the multi indices $(j,\ell,k)$ and
$(i,m,n)$ with $P$ and $Q$, respectively. This way we write
$\psi_P=\psi_{j,\ell,k}$ or $\psi_Q=\psi_{i,m,n}$. Also, we let
$x_P$ and $x_Q$ be the lower left corners $A^{-j}B^{-\ell}k$ and
$A^{-i}B^{-m}n$ of the ``cubes" $P=Q_{j,\ell,k}$ and $Q=Q_{i,m,n}$,
respectively. Let $B_r(x)$ be the Euclidean ball centered in $x$
with radius $r$.

The elements of the affine collection
$$\mathcal{A}_{AB}:=\{\psi_{j,\ell,k}(x)= \abs{\text{det }A}^{j/2}\psi(B^\ell A^j x-k):
    j\geq 0, -2^j\leq\ell\leq 2^j, k\in\mathbb{Z}^2\},$$
have Fourier transform
$$(\psi_{j,\ell,k})^\wedge(\xi)= \abs{\text{det } A}^{-j/2} \hat{\psi}(\xi A^{-j}B^{-\ell})
    \mathbf{e}^{-2\pi i \xi A^{-j}B^{-\ell} k}.$$
Using the anisotropic dilation it is also easy to verify that
$$\psi_{A^{-j}B^{-\ell}}(x-A^{-j}B^{-\ell}k)=\abs{\text{det }A}^{j/2}\psi_{j,\ell,k}(x)
    =\abs{P}^{-1/2}\psi_P(x)$$
and thus
$$\left(\psi_{A^{-j}B^{-\ell}}(\cdot - A^{-j}B^{-\ell}k)\right)^\wedge(\xi)
    =\hat{\psi}(\xi A^{-j}B^{-\ell})\mathbf{e}^{-2\pi i\xi A^{-j}B^{-\ell}k}.$$
We also have
\begin{eqnarray}\label{e:InnrProd_equiv_Conv_Shrlts}
\nonumber
  \ip{f}{\psi_P}
    &=& \ip{f}{\psi_{j,\ell,k}} \\
\nonumber
    &=& \int_{\mathbb{R}^2} f(x) \overline{2^{-3j/2}\psi_{A^{-j}B^{-\ell}}(x-A^{-j}B^{-\ell}k)} dx \\
    &=& \abs{P}^{1/2}(f\ast\tilde{\psi}_{A^{-j}B^{-\ell}})(x_P).
\end{eqnarray}

\vskip0.5cm
\subsection{Almost Orthogonality}\label{sS:Notatn_&_Almost-Orthgnlty}
From the support condition on $\hat{\psi}_1$, the definition of the
matrix $A$ and (\ref{e:Discrt_Shrlt_Cond_Cone_1}), the set of all
shearlets at scale $j$ (for all shear and translation parameters)
interacts with the sets of all shearlets only at scales $j-1,j$ and
$j+1$ (for all shear and translation parameters). The next result
(for functions in $\mathcal{S}$ not necessarily shearlets) is proved
in Subsection \ref{sS:Proofs_Notation}.
\begin{lem}\label{l:c:l:conv_shrlts}
Let $g,h\in \mathcal{S}$. For $i=j-1,\;j,\;j+1\geq 0$, let $Q$ be
identified with $(i,m,n)$. Then, for every $N>2$, there exists a
$C_N>0$ such that
$$\abs{g_{A^{-j}B^{-\ell}}\ast h_Q(x)}\leq
\frac{C_N \abs{Q}^{-\frac{1}{2}}}{(1+2^i\abs{x-x_Q})^{N}},$$ for all
$x\in\mathbb{R}^2$.
\end{lem}

By construction, for the specific case of the ``shearlets on the
cone" we even have the next more informative property stated in the
Fourier domain. The next result is also proved in Subsection
\ref{sS:Proofs_Notation}.
\begin{lem}\label{l:Overlap_bnd}
Let $\text{supp }\hat{\psi}$ be as in Subsection
\ref{sS:Discrt_Shrlts_Cone}. Then, the support of a horizontal
$(\psi_{j,\ell,k})^\wedge$ overlaps with the support of at most 11
other horizontal $(\psi_{i,m,n})^\wedge$ for $(j,\ell)\not =(i,m)$
and all $k,n\in \mathbb{Z}^2$.
\end{lem}

\begin{rem}\label{r:l:Overlap_bnd1}
Since the translation parameters $k$ and $n$ do not affect the
support in the frequency domain, then for
$$f=T_\psi\mathbf{s}=\sum_{Q\in \mathcal{Q}_{AB}}s_Q\psi_Q
=\sum_{i\geq
0}\sum_{m=-2^i,\ldots,2^i}\sum_{n\in\mathbb{Z}^2}s_{i,m,n}\psi_{i,m,n},$$
we formally have that
$$(\tilde{\psi}_{A^{-j}B^{-\ell}}\ast f)(x)=\sum_{i=j-1}^{j+1}\sum_{m(\ell,i)}\sum_{Q\in\mathcal{Q}^{i,m}}
     s_Q(\tilde{\psi}_{A^{-j}B^{-\ell}}\ast \psi_Q)(x),$$
where $m(\ell,i)$ are the shear indices of those shearlets in the
Fourier domain ``surrounding" the support of
$(\tilde{\psi}_{A^{-j}B^{-\ell}})^\wedge$ and the sum
$\sum_{i=j-1}^{j+1}\sum_{m(\ell,i)}$ has at most 11+1 terms for all
$j$ by Lemma \ref{l:Overlap_bnd}.
\end{rem}

\begin{rem}\label{r:l:Overlap_bnd2}
From Lemma \ref{l:Overlap_bnd} the number of horizontal/vertical
shearlets overlapping on the Fourier domain is bounded for all
scales, since the vertical system for $\mathcal{D}^v$ is an
orthonormal rotation of the horizontal system for $\mathcal{D}^h$,
leaving all distances and angles of the supports unaltered.
\end{rem}

\vskip1cm
\section{The characterization}\label{S:T-L-Shearlets}
After defining the distribution spaces we will work on, we will
ignore the ``directions" of the horizontal and vertical
($\mathfrak{d}=\{h,v\}$) shearlets as done in the wavelets case. We
will also ignore the coarse function $\varphi$ and associated
sequence since they are already treated in the literature (see
Section 12 in \cite{FJ90}).

\vskip0.5cm
\subsection{$AB$-anisotropic inhomogeneous Triebel-Lizorkin spaces}\label{sS:Defs_T-L_fun&seq_spcs}
Let $\varphi, \psi\in \mathcal{S}$ be as in Subsection
\ref{sS:Discrt_Shrlts_Cone}.

\begin{defi}\label{d:T-L_func_spcs}
Let $\alpha\in\mathbb{R}$, $0<p<\infty$ and $0<q\leq \infty$. The
\textbf{$AB$-anisotropic inhomogeneous} Triebel-Lizorkin
\textbf{distribution} spaces $\mathbf{F}^{\alpha,q}_p(AB)$ are
defined as the collection of all $f\in\mathcal{S}'$ such that
\begin{eqnarray}\label{e:def_T-L_func_spcs}
\nonumber
  \norm{f}_{\mathbf{F}^{\alpha,q}_p(AB)}
    &=& \norm{f\ast \varphi}_{L^p} \\
    &+& \norm{\left(\sum_{\mathfrak{d}=\{h,v\}}\left\{\sum_{j\geq 0}\sum_{\ell=-2^j}^{2^j}
        [\abs{Q_j}^{-\alpha}\abs{\tilde{\psi}^{\mathfrak{d}}_{A^{-j}_{\mathfrak{d}}B^{-\ell}_{\mathfrak{d}}}\ast
        f}]^q\right\}\right)^{1/q}}_{L^p}<\infty.
\end{eqnarray}

\end{defi}

To work in the sequence level with the shearlets coefficients we
also have the next definition.

\begin{defi}\label{d:T-L_seq_spcs}
Let $\alpha\in\mathbb{R}$, $0<p<\infty$ and $0<q\leq \infty$.  The
\textbf{$AB$-anisotropic inhomogeneous} Triebel-Lizorkin
\textbf{sequence} spaces $\mathbf{f}^{\alpha,q}_p(AB)$ are defined
as the collection of all complex-valued sequences
$\mathbf{s}=\{s_Q\}_{Q\in\mathcal{Q}_{AB}}$ such that
\begin{equation}\label{e:def_T-L_seq_spcs}
  \norm{\mathbf{s}}_{\mathbf{f}^{\alpha,q}_p(AB)}
    = \norm{\left(\sum_{Q\in\mathcal{Q}_{AB}}(\abs{Q}^{-\alpha}\abs{s_Q}\tilde{\chi}_Q)^q\right)^{1/q}}_{L^p}<\infty.
\end{equation}
\end{defi}

We also formally define the \textbf{analysis} and \textbf{synthesis
operators} as
\begin{equation}\label{e:S-T_Ops}
S_\psi f=\{\ip{f}{\psi_Q}\}_{Q\in\mathcal{Q}_{AB}}\;\;\; \text{and }
\;\;\;T_\psi \mathbf{s}=\sum_{Q\in\mathcal{Q}_{AB}}s_Q \psi_Q,
\end{equation}
respectively.

\begin{rem}\label{r:PrsvlFrm_Loss}
Observe that in (\ref{e:def_T-L_func_spcs}) there are no trace of
the characteristic functions $\chi_{\mathcal{D}^h}$,
$\chi_{\mathcal{D}^v}$ and $\chi_{\mathcal{R}}$ (in the Fourier
domain) which enable the identity in $L^2(\mathbb{R}^2)$ via
(\ref{e:Shrlt_Repr_Sys}). Instead, we will simply bound the
operators, since ignoring $\chi_{\mathcal{D}^h}$,
$\chi_{\mathcal{D}^v}$ and $\chi_{\mathcal{R}}$ in
(\ref{e:Shrlt_Repr_Sys}) affects only the Parseval condition on the
frame (see Lemma \ref{l:Overlap_bnd} and Remark
\ref{r:l:Overlap_bnd2}).
\end{rem}

\vskip0.5cm
\subsection{Two basic results}\label{sS:Two_bsc_reslts}
As aforementioned, for the proof of our main result (Theorem
\ref{t:S-T-psi-ops_Bnd}) we follow \cite{FJ90}. This is based on a
kind of Peetre's inequality to bound
$S_\psi:\mathbf{F}^{\alpha,q}_p(AB)\rightarrow
\mathbf{f}^{\alpha,q}_p(AB)$, and a characterization of
$\textbf{f}_p^{\alpha,q}(AB)$ to bound $T_\psi:
\mathbf{f}^{\alpha,q}_p(AB)\rightarrow \mathbf{F}^{\alpha,q}_p(AB)$.
We start with a definition and a well known result.

\begin{defi}\label{d:H-L_Max_Fcn}
The \textbf{Hardy-Littlewood maximal function}, $\mathcal{M}f(x)$,
is given by
$$\mathcal{M}f(x)=\sup_{r>0}\frac{1}{\abs{B_r(x)}}\int_{B_r(x)} \abs{f(y)} dy,$$
for a locally integrable function $f$ on $\mathbb{R}^2$ and where
$B_r(x)$ is the ball with center in $x$ and radius $r$.
\end{defi}
\noindent It is well known that $\mathcal{M}$ is bounded on $L^p$,
$1<p\leq \infty$. It is also true that the next vector-valued
inequality holds (see \cite{FeSt71}).
\begin{thm}\label{t:Feff-Stn_Ineq} \textbf{[Fefferman-Stein]}
For $1<p<\infty$ and $1<q\leq \infty$, there exists a constant
$C_{p,q}$ such that
$$\norm{\left\{\sum_{i=1}^\infty (\mathcal{M}f_i)^q\right\}^{1/q}}_{L^p}
\leq C_{p,q}\norm{\left\{\sum_{i=1}^\infty
f_i^q\right\}^{1/q}}_{L^p},$$ for any sequence $\{f_i:
i=1,2,\ldots\}$ of locally integrable functions.
\end{thm}

\noindent Let us define the \textbf{\emph{shear anisotropic}}
\textbf{Peetre's maximal function}. For all $\lambda>0$,
\begin{equation}\label{e:def_sup_of_conv}
(\psi^{\ast\ast}_{j,\ell,\lambda}f)(x):=\sup_{y\in\mathbb{R}^2}
\frac{\abs{(\psi_{A^{-j}B^{-\ell}}\ast f)(x-y)}}{(1+\abs{B^\ell
A^jy})^{2\lambda}}.
\end{equation}
We then have a \textbf{\emph{shear anisotropic }}\textbf{Peetre's
inequality}. Next lemma is proved in Subsection
\ref{sS:Proof_Two_bsc_reslts}.

\begin{lem}\label{l:SupConv_leq_HLMaxFcn} Let $\psi$ be band
limited and $f\in \mathcal{S}'$. Then, for any real $\lambda>0$,
there exists a constant $C_\lambda$ such that
$$(\psi^{\ast\ast}_{j,\ell,\lambda}f)(x)
\leq C_\lambda\left\{\mathcal{M}(\abs{\psi_{A^{-j}B^{-\ell}}\ast
f}^{1/\lambda})(x)\right\}^\lambda, \;\;\; x\in\mathbb{R}^2.$$
\end{lem}

Identify $Q$ and $P$ with $(i,m,n)$ and $(j,\ell,k)$, respectively.
For all $r>0$, $N\in\mathbb{N}$ and $i\geq j\geq 0$, define
$$(s_{r,N}^\ast)_Q:=
\left(\sum_{P\in\mathcal{Q}^{j,\ell}}\frac{\abs{s_P}^r}{(1+2^j\abs{x_Q-x_P})^N}\right)^{1/r},$$
and
$\mathbf{s}^\ast_{r,N}=\{(s^\ast_{r,N})_Q\}_{Q\in\mathcal{Q}_{AB}}$.
We then have the characterization of the sequence spaces
$\mathbf{f}^{\alpha,q}_p(AB)$ in terms of $\mathbf{s}^\ast_{r,N}$
which is used to prove the boundedness of
$T_\psi:\mathbf{f}^{\alpha,q}_p(AB)\rightarrow
\mathbf{F}^{\alpha,q}_p(AB)$. Next result is also proved in
Subsection \ref{sS:Proof_Two_bsc_reslts}.
\begin{lem}\label{l:s_ast_bnd_s}
Let $\alpha\in \mathbb{R}$, $0<p<\infty$ and $0<q\leq\infty$. Then,
for all $r>0$ and $N>3\max(1,r/q,r/p)$ there exists $C>0$ such that
$$\norm{\mathbf{s}}_{\mathbf{f}^{\alpha,q}_p(AB)}\leq\norm{\mathbf{s}^\ast_{r,N}}_{\mathbf{f}^{\alpha,q}_p(AB)}
\leq C \norm{\mathbf{s}}_{\mathbf{f}^{\alpha,q}_p(AB)}.$$
\end{lem}

\vskip0.5cm
\subsection{Boundedness of $S_\psi$ and $T_\psi$}\label{sS:S-T_psi-ops_Bnd}
As previously mentioned, since we are leaving aside the
characteristic functions $\chi_{\mathcal{D}^h}$,
$\chi_{\mathcal{D}^v}$ and $\chi_{\mathcal{R}}$ of
(\ref{e:Shrlt_Repr_Sys}) one cannot hope for a reproducing identity
for the spaces $\mathbf{F}^{\alpha,q}_p(AB)$.


\begin{thm}\label{t:S-T-psi-ops_Bnd} Let $\alpha\in\mathbb{R}$,
$0<p<\infty$ and $0<q\leq\infty$. Then, the operators $S_\psi :
\mathbf{F}^{\alpha,q}_p(AB)\rightarrow \mathbf{f}^{\alpha,q}_p(AB)$
and $T_\psi : \mathbf{f}^{\alpha,q}_p(AB)\rightarrow
\mathbf{F}^{\alpha,q}_p(AB)$ are well defined and bounded.


\end{thm}
\textbf{Proof.} We prove only the case $q<\infty$. To prove the
boundedness of $S_\psi$ suppose $f\in \mathbf{F}^{\alpha,q}_p(AB)$.
Let $P$ be identified with $(j,\ell,k)$. Then,
$\abs{\tilde{\psi}_{A^{-j}B^{-\ell}}\ast
f(x_P)}\chi_P=\abs{\ip{f}{\psi_P}}\tilde{\chi}_P$, as in
(\ref{e:InnrProd_equiv_Conv_Shrlts}). Let $E=\cup_{\kappa\in
K}Q_{j,\ell,\kappa}$ where $K=\{(0,0),(-1,0),(0,-1),(-1,-1)\}$.
Since $\mathcal{Q}^{j,\ell}$ is a partition of $\mathbb{R}^2$ we
have for $x\in P'$ and $P'\in\mathcal{Q}^{j,\ell}$,
\begin{eqnarray*}
\!\!\!\!\!\!\!\!\!\!\!\!\!\!\!\!\!\!\!\!\!\!\!\!\!\!\!\!\!\!\!\!\!\!\!\!
    & & \!\!\!\!\!\!\!\!\!\!\!\!\!\!\!\!\!\!\!\!\!\!\!\!\!\!\!\!\!\!\!\!\!\!\!\!
        \sum_{P\in\mathcal{Q}^{j,\ell}} [\abs{P}^{-\alpha}\abs{(S_\psi
  f)_P}\tilde{\chi}_P(x)]^q \\
    &=& \abs{\text{det }A}^{j\alpha q} \sum_{P\in\mathcal{Q}^{j,\ell}}
        \left[\abs{\tilde{\psi}_{A^{-j}B^{-\ell}}\ast f(x_P)}\chi_P(x)\right]^q\\
    &\leq& \abs{ \text{det }A}^{j\alpha q} \sum_{P\in\mathcal{Q}^{j,\ell}}
        \sup_{y\in P}\abs{\tilde{\psi}_{A^{-j}B^{-\ell}}\ast
        f(y)}^q \chi_P(x)\\
    &\leq& \abs{ \text{det }A}^{j\alpha q} \sup_{z\in E}
        \abs{\tilde{\psi}_{A^{-j}B^{-\ell}}\ast f(x-z)}^q \\
    &=& \abs{ \text{det }A}^{j\alpha q} \sup_{z\in E}
        \left[\frac{\abs{\tilde{\psi}_{A^{-j}B^{-\ell}}\ast
        f(x-z)}}{(1+\abs{B^\ell A^j z})^{2/\lambda}}\right]^q (1+\abs{B^\ell A^j z})^{q2/\lambda} \\
    &\leq& \abs{\text{det }A}^{j\alpha q}
        \left[\sup_{z\in \mathbb{R}^2} \frac{\abs{\tilde{\psi}_{A^{-j}B^{-\ell}}\ast f(x-z)}}{(1+\abs{B^\ell A^j z})^{2/\lambda}}\right]^q
        \sup_{\kappa\in K}(1+\text{Diam}(Q_{0,0,\kappa}))^{2q/\lambda} \\
    &=& C_{q,\lambda}\abs{\text{det }A}^{j\alpha q}(\tilde{\psi}^{\ast\ast}_{j,\ell,1/\lambda} f)^q(x) \\
    &\leq& C_{q,\lambda}\abs{\text{det }A}^{j\alpha q} \left\{\mathcal{M}\left(\abs{\tilde{\psi}_{A^{-j}B^{-\ell}}\ast
        f}^\lambda\right)(x)\right\}^{q/\lambda},
\end{eqnarray*}
because of Lemma \ref{l:SupConv_leq_HLMaxFcn} (with $1/\lambda$
instead of $\lambda$ in the last inequality). Now, take
$0<\lambda<\min(p,q)$. Then, the previous estimate and Theorem
\ref{t:Feff-Stn_Ineq} yield
\begin{eqnarray*}
  \norm{S_\psi f}_{\mathbf{f}^{\alpha,q}_p(AB)}
    &=& \norm{\left(\sum_{j\geq 0}\sum_{\ell=-2^j}^{2^j}\sum_{P\in \mathcal{Q}^{j,\ell}}
        [\abs{P}^{-\alpha}\abs{(S_\psi f)_P}\tilde{\chi}_P]^q\right)^{1/q}}_{L^p} \\
    &\leq& C\norm{\left(\sum_{j\geq 0}\sum_{\ell=-2^j}^{2^j}
        \left\{\mathcal{M}\left(\abs{\text{det }A}^{j\alpha\lambda}\abs{\tilde{\psi}_{A^{-j}B^{-\ell}}\ast
        f}^\lambda
        \right)\right\}^{q/\lambda}\right)^{1/q}}_{L^p} \\
    &=& C\norm{\left(\sum_{j\geq 0}\sum_{\ell=-2^j}^{2^j}
        \left\{\mathcal{M}\left(\abs{\text{det }A}^{j\alpha\lambda}\abs{\tilde{\psi}_{A^{-j}B^{-\ell}}\ast
        f}^\lambda
        \right)\right\}^{q/\lambda}\right)^{\lambda/q}}_{L^{p/\lambda}}^{1/\lambda} \\
    &\leq& C\norm{\left(\sum_{j\geq 0}\sum_{\ell=-2^j}^{2^j}
        \abs{\text{det }A}^{j\alpha q}\abs{\tilde{\psi}_{A^{-j}B^{-\ell}}\ast f}^q
        \right)^{\lambda/q}}_{L^{p/\lambda}}^{1/\lambda} \\
    &=& C\norm{\left(\sum_{j\geq 0}\sum_{\ell=-2^j}^{2^j}
        [\abs{\text{det }A}^{j\alpha}\abs{\tilde{\psi}_{A^{-j}B^{-\ell}}\ast
        f}]^q
        \right)^{1/q}}_{L^{p}} \\
    &=& C\norm{f}_{\mathbf{F}^{\alpha,q}_p(AB)}.
\end{eqnarray*}

To prove the boundedness of $T_\psi$ suppose
$\mathbf{s}=\{s_Q\}_Q\in \mathbf{f}^{\alpha,q}_p$ and
$f=T_\psi\mathbf{s}=\sum_{Q\in \mathcal{Q}_{AB}}s_Q\psi_Q$. By Lemma
\ref{l:Overlap_bnd} (see also Remark \ref{r:l:Overlap_bnd1}) and
Lemma \ref{l:c:l:conv_shrlts}, we have for $x\in Q'$ and
$Q'\in\mathcal{Q}^{i,m}$,
\begin{eqnarray*}
  \abs{\tilde{\psi}_{A^{-j}B^{-\ell}}\ast f(x)}
    &\leq& \sum_{i=j-1}^{j+1}\sum_{m(\ell,i)}\sum_{Q\in\mathcal{Q}^{i,m}}
        \abs{s_Q}\abs{\tilde{\psi}_{A^{-j}B^{-\ell}}\ast \psi_Q(x)} \\
    &\leq& C\sum_{i=j-1}^{j+1}\sum_{m(\ell,i)}\sum_{Q\in\mathcal{Q}^{i,m}}
        \abs{s_Q}\frac{\abs{Q}^{-1/2}}{(1+2^i\abs{x-x_Q})^N} \\
    &\leq& C'\sum_{i=j-1}^{j+1}\sum_{m(\ell,i)}\sum_{Q\in\mathcal{Q}^{i,m}}
        \abs{s_Q}\frac{\abs{Q}^{-1/2}}{(1+2^i\abs{x_{Q'}-x_Q})^N} \\
    &=& C'\sum_{i=j-1}^{j+1}\sum_{m(\ell,i)}
        \abs{Q}^{-1/2}(s_{1,N}^\ast)_{Q'}\chi_{Q'}(x)\\
    &=& C'\sum_{i=j-1}^{j+1}\sum_{m(\ell,i)}\sum_{Q\in\mathcal{Q}^{i,m}}
        (s_{1,N}^\ast)_Q \tilde{\chi}_Q(x),
\end{eqnarray*}
for all $N>2$ and because $\mathcal{Q}^{i,m}$ is a partition of
$\mathbb{R}^2$. Let $N>3\max(1,1/q,1/p)$. Then, since the pair
$(i,m)$ runs over each pair of scale and shear parameters at most 12
times by Lemma \ref{l:Overlap_bnd}, the previous estimate yields

\begin{eqnarray*}
        \norm{T_\psi \mathbf{s}}_{\mathbf{F}^{\alpha,q}_p(AB)}
    &=& \norm{\left(\sum_{j\geq 0}\sum_{\ell=-2^j}^{2^j}
        [\abs{Q_j}^{-\alpha}\abs{\tilde{\psi}_{A^{-j}B^{-\ell}}\ast f}]^q\right)^{1/q}}_{L^p} \\
    &\leq& C \norm{\left(\sum_{j\geq 0}\sum_{\ell=-2^j}^{2^j}
        \left[\abs{Q_j}^{-\alpha}
        \sum_{i=j-1}^{j+1}\sum_{m(\ell,i)}\sum_{Q\in\mathcal{Q}^{i,m}}
        (s_{1,N}^\ast)_Q\tilde{\chi}_Q\right]^q\right)^{1/q}}_{L^p} \\
    &\leq& C\norm{\left(\sum_{j\geq 0}\sum_{\ell=-2^j}^{2^j}
        \left[\sum_{Q\in\mathcal{Q}^{j,\ell}}
        \abs{Q}^{-\alpha}(s_{1,N}^\ast)_Q\tilde{\chi}_Q\right]^q\right)^{1/q}}_{L^p} \\
    &=& C \norm{\left(\sum_{j\geq 0}\sum_{\ell=-2^j}^{2^j}
        \sum_{Q\in\mathcal{Q}^{j,\ell}}[
        \abs{Q}^{-\alpha}(s_{1,N}^\ast)_Q\tilde{\chi}_Q]^q\right)^{1/q}}_{L^p}\\
    &=& C\norm{\mathbf{s}_{1,N}^\ast}_{\mathbf{f}^{\alpha,q}_p(AB)}
        \leq
        C\norm{\mathbf{s}}_{\mathbf{f}^{\alpha,q}_p(AB)},
\end{eqnarray*}
because $\mathcal{Q}^{j,\ell}$ is a partition of $\mathbb{R}^2$ and
Lemma \ref{l:s_ast_bnd_s} in the last inequality.

\hfill $\blacksquare$ \vskip .5cm   


\begin{rem}\label{r:Def_T-LSpcs_Indep_psi}
With the same arguments as in Remark 2.6 in \cite{FJ90}, the
definition of $\mathbf{F}^{\alpha,q}_p(AB)$ is independent of the
choice of $\psi\in\mathcal{S}$ as long as it satisfies the
requirements in Subsection \ref{sS:Discrt_Shrlts_Cone}.
\end{rem}

\vskip1cm
\section{The identity with smooth Parseval frames}\label{S:Ident_SmthPrsvlFrms}
Recently, Guo and Labate in \cite{GL11} found a way to overcome the
use of characteristic functions in the Fourier domain to restrict
the horizontal/vertical shearlets to the respective cone (see
(\ref{e:Shrlt_Repr_Sys})). The use of these characteristic functions
affects the smoothness of the \emph{boundary shearlets} (those with
$\ell=\pm 2^j$). They slightly modify the definition of these
boundary shearlets instead of projecting them into the cone. This
new shearlets system is not affine-like. However, they do produce
the same frequency tiling as that in Subsection
\ref{sS:Discrt_Shrlts_Cone}.

\vskip0.5cm
\subsection{The new smooth shearlets system}\label{sS:SmthPrsvlFrms}
This subsection is a brief summary of some results in \cite{GL11}
and is intended to show the construction of such smooth Parseval
frames. Let $\phi$ be a $C^\infty$ univariate function such that
$0\leq\phi\leq 1$, with $\hat{\phi}=1$ on $[-1/16,1/16]$ and
$\hat{\phi}=0$ outside $[-1/8,1/8]$ (\emph{i.e.}, $\phi$ is a
rescaled Meyer wavelet). For $\xi\in\hat{\mathbb{R}}^2$, let
$\hat{\Phi}(\xi)=\hat{\phi}(\xi_1)\hat{\phi}(\xi_2)$ and
$W^2(\xi)=\hat{\Phi}^2(2^{-2}\xi)-\hat{\Phi}^2(\xi)$. It follows
that
$$\hat{\Phi}(\xi)+\sum_{j\geq 0}W^2(2^{-2j}\xi)=1, \;\;\text{ for all } \xi\in\hat{\mathbb{R}}^2.$$
Let now $v\in C^\infty(\mathbb{R})$ be such that $v(0)=1$,
$v^{(n)}(0)=0$ for all $n\geq 1$, supp $v\subset [-1,1]$ and
$$\abs{v(u-1)}^2+\abs{v(u)}^2+\abs{v(u+1)}^2=1, \;\; \abs{u}\leq 1.$$
Then, for any $j\geq 0$,
$$\sum_{m=-2^j}^{2^j} \abs{v(2^ju-m)}^2=1, \;\;\abs{u}\leq 1.$$
See Subsection \ref{sS:Discrt_Shrlts_Cone} for comments on the
construction of theses functions and similar properties.

With $V_h(\xi_1,\xi_2)=v(\frac{\xi_2}{\xi_1})$,
$\xi\in\mathcal{D}^h$, the horizontal shearlet system for
$L^2(\mathbb{R}^2)$ is defined as the countable collection of
functions
$$\{\psi^h_{j,\ell,k}:j\geq 0,\abs{\ell}<2^j, k\in\mathbb{Z}^2\},$$
whose elements are defined by their Fourier transform
\begin{equation}\label{e:def_new_Shrlt_Sys_Fourier}
(\psi^h_{j,\ell,k})^\wedge(\xi)=\abs{\text{det } A_h}^{-j/2}
    W(2^{-2j}\xi)V_h(\xi A_h^{-j}B_h^{-\ell})\mathbf{e}^{-2\pi i \xi
    A_h^{-j}B_h^{-\ell}k}, \;\;\; \xi\in\mathcal{D}^h,
\end{equation}
where $A_h$ and $B_h$ are as in Subsection
\ref{sS:Discrt_Shrlts_Cone}. Similarly, one can construct the
vertical shearlet system as in Subsection
\ref{sS:Discrt_Shrlts_Cone}.

For the \emph{boundary shearlets} let $j\geq 1$, $\ell=\pm 2^j$ and
$k\in\mathbb{Z}^2$, then $(\psi_{j,\ell,k})^\wedge(\xi)=
2^{-\frac{3}{2}j-\frac{1}{2}}
W(2^{-2j}\xi)v(2^j\frac{\xi_2}{\xi_1}-\ell)\mathbf{e}^{-2\pi i\xi
2^{-1}A_h^{-j}B_h^{-\ell}k}$ for $\xi\in \mathcal{D}^h$, and
$(\psi_{j,\ell,k})^\wedge(\xi)=2^{-\frac{3}{2}j-\frac{1}{2}}$
$W(2^{-2j}\xi)v(2^j\frac{\xi_1}{\xi_2}-\ell)\mathbf{e}^{-2\pi i\xi
2^{-1}A_v^{-j}B_v^{-\ell}k}$ for $\xi\in \mathcal{D}^v$. When $j=0$,
$\ell=\pm 1$ and $k\in\mathbb{Z}^2$ define
$(\psi_{0,\ell,k})^\wedge(\xi)=
    W(\xi)v(\frac{\xi_2}{\xi_1}-\ell)\mathbf{e}^{-2\pi i\xi k}$
for $\xi\in \mathcal{D}^h$, and $(\psi_{0,\ell,k})^\wedge(\xi)=
    W(\xi)v(\frac{\xi_1}{\xi_2}-\ell)\mathbf{e}^{-2\pi i\xi k}$
for $\xi\in \mathcal{D}^v$. These boundary shearlets are also
$C^\infty(\hat{\mathbb{R}}^2)$ (see \cite{GL11}). This new system is
not affine-like since the function $W$ is not shear-invariant.
However, as previously mentioned, they generate the same frequency
tiling.

The new smooth Parseval frame condition is now written as (see
Theorem 2.3 in \cite{GL11})
\begin{equation}\label{e:SmthPrsvlFrm_Prop_Shrlts}
\abs{\hat{\Phi}(\xi)}^2+\sum_{\mathfrak{d}=1}^2\sum_{j\geq
0}\sum_{\abs{\ell}<2^j} \abs{\hat{\psi}^{\mathfrak{d}}(\xi
A_\mathfrak{d}^{-j} B_\mathfrak{d}^{-\ell})}^2 + \sum_{j\geq
0}\sum_{\ell=\pm 2^j} \abs{\hat{\psi}(\xi A^{-j} B^{-\ell})}^2=1,
\end{equation}
for all $\xi\in\hat{\mathbb{R}}^2$ and where $\mathfrak{d}=1,2$
stands for horizontal and vertical directions and we omit the
subindex for the matrices of the boundary shearlets. Notice that now
there do not exist characteristic functions as in
(\ref{e:Shrlt_Repr_Sys}).

\vskip0.5cm
\subsection{The reproducing identity on $\mathcal{S}'$}\label{sS:Repr_Ident_Distrbtns}
Our goal is to show that, with the smooth Parseval frames of
shearlets of Guo and Labate in \cite{GL11}, $T_\psi\circ S_\psi$ is
the identity on $\mathcal{S}'$ and, therefore, on
$\mathbf{F}^{\alpha,q}_p(AB)$. First, we show that any
$f\in\mathcal{S}'$ admits a kind of Littlewood-Paley decomposition
with shear anisotropic dilations, for which we follow \cite{BH05}.
Then, we show the reproducing identity in $\mathcal{S}'$ following
\cite{FJW}. Denote $\hat{\Phi}=(\psi_{-1})^\wedge$.

\begin{lem}\label{l:Ident_Conv}
Let $\{\psi_{j,\ell,k}:j\geq 0, \ell=-2^j,\ldots,2^j,
k\in\mathbb{Z}^2\}$ be the smooth shearlet system that verifies
(\ref{e:SmthPrsvlFrm_Prop_Shrlts}). Then, for any
$f\in\mathcal{S}'$,
\begin{eqnarray*}
  f=f\ast\tilde{\psi}_{-1}\ast\psi_{-1}
    &+& \sum_{\mathfrak{d}=1}^2\sum_{j\geq 0}\sum_{\abs{\ell}<2^j}
        f\ast\tilde{\psi}^\mathfrak{d}_{A_\mathfrak{d}^{-j}B_\mathfrak{d}^{-\ell}}\ast\psi^\mathfrak{d}_{A_\mathfrak{d}^{-j}B_\mathfrak{d}^{-\ell}} \\
    &+& \sum_{j\geq 0}\sum_{\ell=\pm 2^j}
        f\ast\tilde{\psi}_{A^{-j}B^{-\ell}}\ast\psi_{A^{-j}B^{-\ell}},
\end{eqnarray*}
with convergence in $\mathcal{S}'$.
\end{lem}
\textbf{Proof.} One can see Peetre's discussion on pp. 52-54 of
\cite{Pe76} regarding convergence. Since the Fourier transform
$\mathcal{F}$ is an isomorphism of $\mathcal{S}'$, it suffices to
show that
\begin{eqnarray*}
  \hat{f}(\xi)=\hat{f}(\xi) \abs{(\psi_{-1})^\wedge(\xi)}^2
    &+& \sum_{\mathfrak{d}=1}^2\sum_{j\geq 0}\sum_{\abs{\ell}<2^j}
        \hat{f}(\xi) \abs{\hat{\psi}^\mathfrak{d}(\xi A_\mathfrak{d}^{-j} B_\mathfrak{d}^{-\ell})}^2\\
    &+& \sum_{j\geq 0}\sum_{\ell=\pm 2^j}
        \hat{f}(\xi) \abs{\hat{\psi}(\xi A^{-j} B^{-\ell})}^2
\end{eqnarray*}
converges in $\mathcal{S}'$. Since the equality is a straight
consequence of (\ref{e:SmthPrsvlFrm_Prop_Shrlts}), we will only show
convergence in $\mathcal{S}'$ of the right-hand side of the equality
for those shearlets with $j\geq 0$ ($\psi_{-1}$ is in fact a scaling
function of a Meyer wavelet). Suppose that $\hat{f}$ has order $\leq
m$. This is, there exists an integer $n\geq 0$ and a constant $C$
such that
$$\abs{\ip{\hat{f}}{g}}\leq C \sup_{\abs{\alpha}\leq n, \abs{\beta}\leq m} \norm{g}_{\alpha,\beta}, \;\;\; \text{ for all } g\in \mathcal{S},$$
where $\norm{g}_{\alpha,\beta}=
\sup_{\xi\in\hat{\mathbb{R}}^2}\abs{\xi^\alpha}\abs{\partial^\beta
g(\xi)}$ denotes the usual semi-norm in $\mathcal{S}$ for
multi-indices $\alpha$ and $\beta$. Then,
$$\abs{\ip{\hat{f}\abs{(\psi_{A^{-j}B^{-\ell}})^\wedge}^2}{g}}=\abs{\ip{\hat{f}}{\abs{(\psi_{A^{-j}B^{-\ell}})^\wedge}^2 g}}
    \leq C \sup_{\abs{\alpha}\leq n, \abs{\beta}\leq m} \norm{\abs{(\psi_{A^{-j}B^{-\ell}})^\wedge}^2 g}_{\alpha,\beta}.$$
As in Lemma 2.5 in \cite{GL07}, one can prove that
$$\sup_{\abs{\beta}=m} \norm{\partial^\beta \abs{(\psi_{A^{-j}B^{-\ell}})^\wedge}^2}_\infty \leq C 2^{-jm}.$$
Hence, by the compact support conditions of
$(\psi_{A^{-j}B^{-\ell}})^\wedge(\xi)$ (see Subsection
\ref{sS:Discrt_Shrlts_Cone})
\begin{eqnarray*}
    & & \!\!\!\!\!\!\!\!\!\!\!\!\!\!\!\!\!\!\!\!\!\!\!\!\!\!\!\!\!\!\!\!\!\!\!\!\!\!\!\!
            \sup_{\abs{\alpha}\leq n, \abs{\beta}\leq m} \norm{\abs{(\psi_{A^{-j}B^{-\ell}})^\wedge}^2
            g}_{\alpha,\beta}\\
    &\leq& C \sup_{\xi\in \hat{\mathbb{R}}^2} \left[(1+\abs{\xi})^n
            \sup_{\abs{\beta}\leq m} \abs{\partial^\beta\abs{(\psi_{A^{-j}B^{-\ell}})^\wedge(\xi)}^2}
            \sup_{\abs{\beta}\leq m} \abs{\partial^\beta g(\xi)}\right] \\
    &\leq& C \sup_{\xi\in \text{supp}(\psi_{A^{-j}B^{-\ell}})^\wedge(\xi)} (1+\abs{\xi})^n \sup_{\abs{\beta}\leq m}\abs{\partial^\beta g(\xi)} \\
    &\leq& C \sup_{\abs{\alpha}\leq n+1, \abs{\beta}\leq m} \norm{g}_{\alpha,\beta}
            \sup_{\xi\in\text{supp}(\psi_{A^{-j}B^{-\ell}})^\wedge(\xi)} (1+\abs{\xi})^{-1}\\
    &\leq& C \sup_{\abs{\alpha}\leq n+1, \abs{\beta}\leq m} \norm{g}_{\alpha,\beta}
            (1+2^{2j-4})^{-1}\leq C2^{-2j},
\end{eqnarray*}
which proves the convergence in $\mathcal{S}'$.

\hfill $\blacksquare$ \vskip .5cm   

\begin{lem}\label{l:equiv_Conv_CctSuppFourier} Let $g\in\mathcal{S}'$ and $h\in\mathcal{S}$ be such that
$$\text{\emph{supp} }\hat{g}, \text{\emph{ supp} }\hat{h}\subset [-1/2,1/2]^2B^\ell A^j=QB^\ell
A^j,
    \;\;\;\;\;\;  j\geq 0, \ell=-2^j,\ldots, 2^j.$$
Then,
$$g\ast h= \sum_{k\in\mathbb{Z}^2} \abs{\text{\emph{det} }A}^{-j} g(A^{-j}B^{-\ell}k) h(x-A^{-j}B^{-\ell}k),$$
with convergence in $\mathcal{S}'$.
\end{lem}
\textbf{Proof.} Suppose first that $g\in\mathcal{S}$. We can express
$\hat{g}$ by its Fourier series as
$$\hat{g}(\xi)=\sum_{k\in\mathbb{Z}^2} \abs{\text{det }A}^{-j/2} \mathbf{e}^{-2\pi i\xi A^{-j}B^{-\ell}k} \cdot
    \left(\int_{QB^\ell A^j} \hat{g}(\omega)\cdot\abs{\text{det }A}^{-j/2}\mathbf{e}^{2\pi i\omega A^{-j}B^{-\ell}k} d\omega\right).$$
By the Fourier inversion formula in $\hat{\mathbb{R}}^2$ we have
$$\hat{g}(\xi)=\sum_{k\in\mathbb{Z}^2} \abs{\text{det }A}^{-j/2} \mathbf{e}^{-2\pi i\xi A^{-j}B^{-\ell}k} \cdot
    g(A^{-j}B^{-\ell}k), \;\;\;\;\;\; \xi\in QB^\ell A^j.$$
Since $\hat{g}$ has compact support, $g(A^{-j}B^{-\ell}k)$ makes
sense (by the Paley-Wiener theorem). Since $\text{supp }
\hat{h}\subset QB^\ell A^j$ and $g\ast h=(\hat{g}\hat{h})^\vee$,
\begin{eqnarray*}
  g\ast h
    &=& \sum_{k\in \mathbb{Z}^2}\abs{\text{det }A}^{-j} g(A^{-j}B^{-\ell}k) [\mathbf{e}^{-2\pi i\xi A^{-j}B^{-\ell}k} \hat{h}(\cdot)]^\vee \\
    &=& \sum_{k\in \mathbb{Z}^2}\abs{\text{det }A}^{-j} g(A^{-j}B^{-\ell}k)
            h(x-A^{-j}B^{-\ell}k),
\end{eqnarray*}
which proves the convergence for $g\in\mathcal{S}$. To remove this
assumption one uses the same standard regularization argument as in
the proof of Lemma \ref{l:Max_Fcn_leq_HL-Max-Fcn}. This
regularization argument is the same used in Lemma(6.10) in
\cite{FJW}.

\hfill $\blacksquare$ \vskip .5cm   

\begin{thm}\label{t:Ident_in_Distributions}
Let the shearlet system $\{\psi_{j,\ell,k}\}$ be constructed as in
Subsection \ref{sS:SmthPrsvlFrms} such that it is a smooth Parseval
frame that verifies (\ref{e:SmthPrsvlFrm_Prop_Shrlts}). The
composition of the analysis and synthesis operators $T_\psi\circ
S_\psi$ (see (\ref{e:S-T_Ops}) for the definitions) is the identity
$$f=\sum_{Q\in\mathcal{Q}_{AB}} \ip{f}{\psi_Q}\psi_Q,$$
in $\mathcal{S}'$.

\end{thm}
\textbf{Proof.} As in (\ref{e:InnrProd_equiv_Conv_Shrlts}),
$f\ast\tilde{\psi}_{A^{-j}B^{-\ell}}(A^{-j}B^{-\ell}k)=
f\ast\tilde{\psi}_{A^{-j}B^{-\ell}}(x_P) = \abs{\text{det
}A}^{j/2}\ip{f}{\psi_P}$, where $P$ is identified with $(j,\ell,k)$.
Also, as in Subsection \ref{sS:Notatn_&_Almost-Orthgnlty},
$\psi_{A^{-j}B^{-\ell}}(x-A^{-j}B^{-\ell}k)=\abs{\text{det
}A}^{j/2}\psi_P(x)$.
Let $g=f\ast\tilde{\psi}_{A^{-j}B^{-\ell}}$ and
$h=\psi_{A^{-j}B^{-\ell}}$. By construction, 
supp$(\psi_{j,\ell,k})^\wedge(\xi)\subset QB^\ell A^j$. Therefore,
Lemma \ref{l:equiv_Conv_CctSuppFourier} yields
\begin{eqnarray*}
  f\ast\tilde{\psi}_{A^{-j}B^{-\ell}}\ast\psi_{A^{-j}B^{-\ell}}
    &=& \sum_{k\in\mathbb{Z}^2} \ip{f}{\psi_{j,\ell,k}}\psi_{j,\ell,k} \\
    &=& \sum_{P\in\mathcal{Q}^{j,\ell}} \ip{f}{\psi_P}\psi_P.
\end{eqnarray*}
By appropriately summing over $\mathfrak{d}=1,2$, $j\geq 0$ and
$\ell=-2^j,\ldots,2^j$, Lemma \ref{l:Ident_Conv} yields the result.

\hfill $\blacksquare$ \vskip .5cm   

\vskip1cm
\section{Relations between $\mathbf{F}^{\alpha_1,q_1}_{p_1}$ and $\mathbf{F}^{\alpha_2,q_2}_{p_2}(AB)$}\label{S:Embeddings}
In this section we prove embeddings of classical dyadic (isotropic)
inhomogeneous Triebel-Lizorkin spaces into the just defined highly
anisotropic inhomogeneous Triebel-Lizorkin spaces, and viceversa,
for certain parameters. We also show that some functions in each of
these spaces vanish in the other spaces for certain parameters. Let
$A$ and $B$ be as in Subsection \ref{sS:Discrt_Shrlts_Cone}. A
dyadic cube will be denoted by $Q$ and a shear anisotropic ``cube"
(a parallelepiped) will be denoted by $P$.

We start with some definitions regarding the classical dyadic spaces
(see Sections 2 and 12 in \cite{FJ90}). Let $\varphi,\theta,
\Phi,\Theta$ be the analyzing and synthesizing functions of the
$\varphi$-transform of Frazier and Jawerth. Then,
$\varphi,\theta,\Phi$ and $\Theta$ satisfy: 1)
$\varphi,\theta,\Phi,\Theta\in\mathcal{S}$, 2) supp $\hat{\varphi}$,
supp
$\hat{\theta}\subset\{\xi\in\hat{\mathbb{R}}^2:\frac{1}{2}\leq\abs{\xi}\leq
2\}$ and supp $\hat{\Phi}$, supp
$\hat{\Theta}\subset\{\xi\in\hat{\mathbb{R}}^2:\abs{\xi}\leq 2\}$,
3) $\abs{\hat{\varphi}(\xi)}, \abs{\hat{\theta}(\xi)}\geq c>0$ if
$\frac{3}{5}\leq\abs{\xi}\leq\frac{5}{3}$ and
$\abs{\hat{\Phi}(\xi)}, \abs{\hat{\Theta}(\xi)}\geq c>0$ if
$\abs{\xi}\leq\frac{5}{3}$, and 4)
$\hat{\tilde{\Phi}}(\xi)\hat{\Theta}(\xi)+\sum_{\nu\in\mathbb{Z}_+}\overline{\hat{\varphi}(2^{-\nu}\xi)}\hat{\theta}(2^{-\nu}\xi)=1$.
Let $\mathcal{D}_+$ denote the set of dyadic cubes with $l(Q)\leq 1$
where $l(Q)$ is the side size of $Q$. Let
$\varphi_{\nu,k}(x)=2^\nu\varphi(2^\nu x-k)$ be the $L^2$-normalized
dilation and $\varphi_{2^\nu I}(x)=2^{2\nu}\varphi(2^\nu x)$, where
$I$ is the identity matrix.

For $\alpha\in\mathbb{R}$, $0<q\leq\infty$, $0<p<\infty$, the
(dyadic) inhomogeneous Triebel-Lizorkin space
$\mathbf{F}^{\alpha,q}_p$ is the collection of all
$f\in\mathcal{S}'$ such that  (see Lemma 12.1 in \cite{FJ90} for a
discussion on the dilation indices)
$$\norm{f}_{\mathbf{F}^{\alpha,q}_p}=\norm{\Phi\ast f}_{L^p}+\norm{\left(\sum_{\mathcal{D}_+} (2^{\nu\alpha}\abs{\varphi_{2^\nu I}\ast f})^q\right)^{1/q}}_{L^p}<\infty.$$

For $\alpha\in\mathbb{R}$, $0<q\leq\infty$, $0<p<\infty$, the
(dyadic) inhomogeneous Triebel-Lizorkin sequence space
$\mathbf{f}^{\alpha,q}_p$ is the collection of all complex-valued
sequences $\mathbf{s}$ such that
$$\norm{\mathbf{s}}_{\mathbf{f}^{\alpha,q}_p}=\norm{\left(\sum_{Q:l(Q)\leq 1} (2^{\nu\alpha}\abs{s_Q}\tilde{\chi}_Q)^q\right)^{1/q}}_{L^p}<\infty,$$
where $\tilde{\chi}_Q(x)=\abs{Q}^{-\frac{1}{2}}\chi_Q(x)$ is the
$L^2$-normalized characteristic function of $Q\in\mathcal{D}_+$.

Let $\mathbf{s}=\{s_Q\}_Q$, where we identify $Q$ with the pair
$(\nu,k)\in\mathbb{Z}_+\times\mathbb{Z}^2$. For $0<r\leq\infty$ and
$\lambda>0$, define the sequence
$\mathbf{s}_{r,\lambda}^\ast=\{(s^\ast_{r,\lambda})_Q\}_{Q\in\mathcal{D}_+}$
by
$$(s_{r,\lambda}^\ast)_{Q'}=\left(\sum_{Q:l(Q)=l(Q')}
    \frac{\abs{s_Q}^r}{(1+l(Q')^{-1}\abs{x_Q-x_{Q'}})^\lambda}\right)^{1/r},$$
where $x_Q=2^{-\nu}k$ is the lower left corner of
$Q_{\nu,k}=2^{-\nu}(Q_0+k)$, see p. 48 of \cite{FJ90} and compare
with the similar definition at Subsection \ref{sS:Two_bsc_reslts}.

\subsection{The embeddings}\label{sS:Embeddings}
We start with a result on almost orthogonality of functions under
highly anisotropic and dyadic dilations.
\begin{lem}\label{l:conv_shrlts-wvlts}
Let $\psi,\varphi\in\mathcal{S}$. For $j\geq 0$, $\abs{\ell}\leq
2^j$ and $k\in\mathbb{Z}^2$,
$$\int_{\mathbb{R}^2} \abs{\psi(B^\ell A^j(x-y))}\abs{\varphi(2^{2j}y)} dy
    \leq \frac{2^{-3j}}{(1+2^j\abs{x})^N},$$
for all $N>2$.
\end{lem}
\textbf{Proof}. Since $\psi,\varphi\in\mathcal{S}$,
\begin{eqnarray*}
    & & \!\!\!\!\!\!\!\!\!\!\!\!\!\!\!\!\!\!\!\!\!\!\!\!\!\!\!\!\!\!\!\!
        \int_{\mathbb{R}^2} \abs{\psi(B^\ell
        A^j(x-y))}\abs{\varphi(2^{2j}y)} dy\\
    &\lesssim& \int_{\mathbb{R}^2}
        \frac{1}{(1+\abs{B^\ell A^j(x-y)})^N}
        \frac{1}{(1+\abs{2^{2j}y})^N}dy.
\end{eqnarray*}
Define
\begin{eqnarray*}
  E_1 &=& \{y\in\mathbb{R}^2:2^j\abs{x-y}\leq 3\} \\
  E_2 &=& \{y\in\mathbb{R}^2:2^j\abs{x-y}> 3, \abs{y}\leq \abs{x}/2\} \\
  E_3 &=& \{y\in\mathbb{R}^2:2^j\abs{x-y}> 3, \abs{y}> \abs{x}/2\}.
\end{eqnarray*}
For $y\in E_1$, $1+2^j\abs{x}\leq 1+2^j\abs{x-y}+2^j\abs{y}\leq
4(1+2^{2j}\abs{y})$. If $y\in E_3$, $1+2^j\abs{x}\leq
1+2^{j+1}\abs{y}\leq 2(1+2^{2j}\abs{y})$. When $y\in E_2$,
$2^{j-1}\abs{x}<2^j(\abs{x}-\abs{y})\leq 2^j\abs{x-y}$, which
implies $4\abs{B^\ell A^j(x-y)}\geq
2^{j-1}\abs{x-y}+3\cdot2^{j-1}\abs{x-y}\geq \frac{3}{2} +
\frac{3}{2}2^{j-1}\abs{x}$ or $8(1+\abs{B^\ell A^j(x-y)})\geq
1+2^j\abs{x}$. Hence,
\begin{eqnarray*}
    & & \!\!\!\!\!\!\!\!\!\!\!\!\!\!\!\!\!\!\!\!\!\!\!\!\!\!\!\!\!\!\!\!
        \int_{\mathbb{R}^2} \abs{\psi(B^\ell
        A^j(x-y))}\abs{\varphi(2^{2j}y)} dy\\
    &\lesssim& \frac{1}{(1+2^j\abs{x})^N} \int_{E_1\cup E_3}
        \frac{1}{(1+\abs{B^\ell A^j(x-y)})^N}dy \\
    & & \;\;\; + \frac{1}{(1+2^j\abs{x})^N} \int_{E_2}
        \frac{1}{(1+2^{2j}\abs{y})^N}dy\\
    &\lesssim& \left[\frac{2^{-3j}}{(1+2^j\abs{x})^N}+\frac{2^{-4j}}{(1+2^j\abs{x})^N}\right]
        \lesssim \frac{2^{-3j}}{(1+2^j\abs{x})^N},
\end{eqnarray*}
for all $N>2$.

\hfill $\blacksquare$ \vskip .5cm   

The definitions of $E_1,E_2,E_3$ in Lema \ref{l:conv_shrlts-wvlts}
allow us to have a ``height" of $2^{-3j}$ and a decreasing of
$(1+2^j\abs{x})^{-N}$. By defining $E_1,E_2,E_3$ as in Lemma
\ref{l:conv_shrlts} would only yield a ``height" of $2^{-2j}$ and a
decreasing of $(1+\abs{x})^{-N}$.

\begin{thm}\label{t:ClssDydc-T-L_in_ShrAnistr-T-L}
Let $\alpha_1,\alpha_2\in \mathbb{R}$, $0<q\leq\infty$, $0<p<\infty$
and $\lambda>2\max(1,1/q,1/p)$. If $3\alpha_2+\frac{1}{q}+\lambda<
\alpha_1$,
$$\mathbf{F}^{\alpha_1,q}_p \hookrightarrow \mathbf{F}^{\alpha_2,q}_p(AB).$$
\end{thm}
\textbf{Proof}. To shorten notation write
$\mathbf{F}_1=\mathbf{F}^{\alpha_1,q}_p$,
$\mathbf{f}_1=\mathbf{f}^{\alpha_1,q}_p$ and
$\mathbf{F}_2=\mathbf{F}^{\alpha_2,q}_p(AB)$. We will actually prove
that, for $f=\sum_{Q\in\mathcal{D}_+}s_Q\varphi_Q\in\mathbf{F}_1$,
$$\norm{f}_{\mathbf{F}_2}\lesssim\norm{\mathbf{s}_{1,\lambda}^\ast}_{\mathbf{f}_1}\lesssim\norm{\mathbf{s}}_{\mathbf{f}_1}\lesssim\norm{f}_{\mathbf{F}_1},$$
where, of course, the inequality we are interested to prove is the
first one and the last two are proved in \cite{FJ90}. From the
compact support conditions of $(\varphi_{\nu,k})^\wedge$ and
$(\psi_{A^{-j}B^{-\ell}})^\wedge$ and their dyadic and highly
anisotropic expansion, respectively, we formally get
$$\psi_{A^{-j}B^{-\ell}}\ast f=
    \sum_{\nu=2j-5}^{2j}\sum_{k\in\mathbb{Z}^2} s_{\nu,k}\psi_{A^{-j}B^{-\ell}}\ast\varphi_{\nu,k}.$$
Therefore, writing $\psi_{A^{-j}B^{-\ell}}(x)=\abs{\text{det
}A}^j\psi(B^\ell A^j x)$ and $\varphi_{\nu,k}(x)=2^\nu\varphi(2^\nu
x-k)$, Lemma \ref{l:conv_shrlts-wvlts} yields
\begin{eqnarray*}
  \norm{f}_{\mathbf{F}_2}
    &=& \norm{\left(\sum_{j\geq 0}2^{3j\alpha_2q}\sum_{\abs{\ell}\leq 2^j}
        [\abs{\sum_{\nu=2j-5}^{2j}\sum_{k\in\mathbb{Z}^2}
        s_{\nu,k}\psi_{A^{-j}B^{-\ell}}\ast\varphi_{\nu,k}(\cdot)}]^q\right)^{1/q}}_{L^p} \\
    &\lesssim& \norm{\left(\sum_{j\geq 0}2^{3j\alpha_2q}\sum_{\abs{\ell}\leq 2^j}
        [\abs{\sum_{k\in\mathbb{Z}^2}
        s_{2j,k}\psi_{A^{-j}B^{-\ell}}\ast\varphi_{2j,k}(\cdot)}]^q\right)^{1/q}}_{L^p} \\
    &\lesssim& \norm{\left(\sum_{j\geq 0}2^{3j\alpha_2q}\sum_{\abs{\ell}\leq 2^j}
        [\sum_{k\in\mathbb{Z}^2}
        \abs{s_{2j,k}}\frac{2^{2j}}{(1+2^j\abs{\cdot+2^{-2j}k})^N}]^q\right)^{1/q}}_{L^p},
\end{eqnarray*}
for all $N>2$. Let $\lambda>2\max(1,r/q,r/p)$ for some $r>0$. 
Following the proof of the second part of Theorem
\ref{t:S-T-psi-ops_Bnd}, if $x\in Q'$ and $Q'\in\mathcal{D}^{2j}$,
\begin{eqnarray*}
    & &\!\!\!\!\!\!\!\!\!\!\!\!\!\!\!\!\!\!\!\!\!\!\!\!\!\!\!\!\!\!\!\!
        \sum_{k\in\mathbb{Z}^2}
        \frac{\abs{s_{2j,k}}2^{2j}}{(1+2^j\abs{x-2^{-2j}k})^\lambda}\\
    &=& \sum_{k\in\mathbb{Z}^2} \frac{\abs{s_{2j,k}}2^{2j}
        \cdot
        2^{j\lambda}}{2^{j\lambda}(1+2^j\abs{x-2^{-2j}k})^\lambda}
        \leq 2^{j\lambda}\sum_{k\in\mathbb{Z}^2}
        \frac{\abs{s_{2j,k}}2^{2j}}{(1+2^{2j}\abs{x-2^{-2j}k})^\lambda}\\
    &\lesssim& 2^{j\lambda}\sum_{Q\in\mathcal{D}^{2j}} \abs{Q}^{-\frac{1}{2}}\abs{(s_{1,\lambda}^\ast)_Q}\chi_Q(x)
        = 2^{j\lambda}\sum_{Q\in\mathcal{D}^{2j}}
        \abs{(s_{1,\lambda}^\ast)_Q}\tilde{\chi}_Q(x),
\end{eqnarray*}
since $\mathcal{D}^{2j}$ is a partition of $\mathbb{R}^2$. Hence,
\begin{eqnarray*}
  \norm{f}_{\mathbf{F}_2}
    &\lesssim& \norm{\left(\sum_{j\geq 0}2^{3j\alpha_2q}(2^{j+1}+1)
        [2^{j\lambda}\sum_{Q\in\mathcal{D}^{2j}}
        \abs{(s_{1,\lambda}^\ast)_Q}\tilde{\chi}_Q(\cdot)]^q\right)^{1/q}}_{L^p} \\
    &\lesssim& \norm{\left(\sum_{j\geq 0}\sum_{Q\in\mathcal{D}^{2j}}
        [2^{3j\alpha_2+\frac{j}{q}+j\lambda}\abs{(s_{1,\lambda}^\ast)_Q}\tilde{\chi}_Q(\cdot)]^q\right)^{1/q}}_{L^p}\\
    &\lesssim& \norm{\left(\sum_{j\geq 0}\sum_{Q\in\mathcal{D}^{j}}
        [2^{3j\alpha_2+\frac{j}{q}+j\lambda}\abs{(s_{1,\lambda}^\ast)_Q}\tilde{\chi}_Q(\cdot)]^q\right)^{1/q}}_{L^p} \\
    &\lesssim& \norm{\left(\sum_{Q\in\mathcal{D}_+}
        [\abs{Q}^{-\frac{\alpha_1}{2}}\abs{(s_{1,\lambda}^\ast)_Q}\tilde{\chi}_Q(\cdot)]^q\right)^{1/q}}_{L^p}
        = \norm{\mathbf{s}^\ast_{1,\lambda}}_{\mathbf{f}_1},
\end{eqnarray*}
since $2\mathbb{N}\subset\mathbb{N}$ and
$2^{3j\alpha_2+\frac{j}{q}+j\lambda} \leq
2^{j\alpha_1}=\abs{Q}^{-\frac{\alpha_1}{2}}$ for a
$Q\in\mathcal{D}^j$. Let now $r=\min(q,p)$. Since
$\lambda/2>\max(1,r/q,r/p)$, $r/(\lambda/2)<\min(r,q,p)$. We can
choose $a$ such that $r/(\lambda/2)<a<\min(r,q,p)$. Then,
$0<a<r<\infty$, $\lambda>2r/a$, $q/a>1$ and $p/a>1$. Following the
proof of Lemma 2.3 in \cite{FJ90} (in which $a$ and $\lambda$ are
defined differently) we get
$\norm{\mathbf{s}_{1,\lambda}^\ast}_{\mathbf{f}_1}
    \lesssim\norm{\mathbf{s}}_{\mathbf{f}_1}$,
and from Theorem 2.2 in \cite{FJ90},
$\norm{\mathbf{s}}_{\mathbf{f}_1}\lesssim\norm{f}_{\mathbf{F}_1}$,
which finishes the proof.
\hfill $\blacksquare$ \vskip .5cm   

\begin{thm}\label{t:ShrAnistr-T-L_in_ClssDydc-T-L}
Let $\alpha_1,\alpha_2\in \mathbb{R}$, $0<q\leq\infty$, and
$0<p<\infty$. If $\alpha_1+1\leq 3\alpha_2$,
$$\mathbf{F}^{\alpha_2,q}_p(AB) \hookrightarrow \mathbf{F}^{\alpha_1,q}_p.$$
\end{thm}
\textbf{Proof.} To shorten notation write
$\mathbf{F}_1=\mathbf{F}^{\alpha_1,q}_p$,
$\mathbf{F}_2=\mathbf{F}^{\alpha_2,q}_p(AB)$ and
$\mathbf{f}_2=\mathbf{f}^{\alpha_2,q}_p$. Suppose
$f=\sum_{P\in\mathcal{Q}_{AB}}s_P\psi_P\in \textbf{F}_2$ and let
$\lambda>3\max(1,1/q,1/p)$. From the compact support conditions on
$(\varphi_{2^\nu I})^\wedge$ and $(\psi_{j,\ell,k})^\wedge$, $j\sim
\lfloor\nu/2\rfloor$. Also, for $\nu\geq 0$, $\lfloor\nu/2\rfloor$
runs through $\mathbb{Z}_+$ twice. So, writing
$\psi_{\nu,\ell,k}(x)=\abs{\text{det }A}^{\nu/2}\psi(B^\ell A^\nu x
-k)$ and $\varphi_{2^{2\nu}I}(x)=2^{4\nu}\varphi(2^{2\nu} x)$,
\begin{eqnarray*}
  \norm{f}_{\textbf{F}_1}
    &=& \norm{\left(\sum_{\nu\geq 0} (\abs{Q_\nu}^{-\alpha_1/2}\abs{\varphi_{2^\nu I}\ast f(\cdot)})^q\right)^{1/q}}_{L^p} \\
    &\lesssim& \norm{\left(\sum_{\nu\geq 0} 2^{\nu\alpha_1
        q}(\abs{\sum_{\abs{\ell}\leq2^{\lfloor\nu/2\rfloor}}\sum_{k\in\mathbb{Z}^2}s_{\lfloor\nu/2\rfloor,\ell,k}\cdot
        \varphi_{2^\nu I}\ast \psi_{\lfloor\nu/2\rfloor,\ell,k}(\cdot)})^q\right)^{1/q}}_{L^p}  \\
    &\lesssim& 2\norm{\left(\sum_{\nu\geq 0} 2^{\nu\alpha_1
        q}(\abs{\sum_{\abs{\ell}\leq2^\nu}\sum_{k\in\mathbb{Z}^2}s_{\nu,\ell,k}\cdot
        \varphi_{2^{2\nu}I}\ast \psi_{\nu,\ell,k}(\cdot)})^q\right)^{1/q}}_{L^p} \\
    &\lesssim& \norm{\left(\sum_{\nu\geq 0} 2^{\nu\alpha_1
        q}(\sum_{\abs{\ell}\leq2^\nu}\sum_{k\in\mathbb{Z}^2}\abs{s_{\nu,\ell,k}}\cdot
        \frac{2^{-3\nu}\cdot2^{4\nu}\cdot2^{3\nu/2}}{(1+2^\nu\abs{\cdot-A^{-\nu}B^{-\ell}k})^N})^q\right)^{1/q}}_{L^p},
\end{eqnarray*}
for all $N>2$, by Lemma \ref{l:conv_shrlts-wvlts}. Continuing as in
the second part of the proof of Theorem \ref{t:S-T-psi-ops_Bnd}, if
$x\in P'$ and $P'\in\mathcal{Q}^{\nu,\ell}$,
\begin{eqnarray*}
  \norm{f}_{\textbf{F}_1}
    &\lesssim& \norm{\left(\sum_{\nu\geq 0} 2^{\nu\alpha_1
        q+\nu q}(\sum_{\abs{\ell}\leq2^\nu}\sum_{P\in\mathcal{Q}^{\nu,\ell}}\abs{s_{P}}\cdot
        \frac{\abs{P}^{-1/2}}{(1+2^\nu\abs{\cdot-x_P})^N})^q\right)^{1/q}}_{L^p} \\
    &\lesssim& \norm{\left(\sum_{\nu\geq 0} 2^{\nu q(\alpha_1 +1)}[\sum_{\abs{\ell}\leq2^\nu}\sum_{P\in\mathcal{Q}^{\nu,\ell}}
        (s_{1,N}^\ast)_P\tilde{\chi}_P(\cdot)]^q\right)^{1/q}}_{L^p},
\end{eqnarray*}
because $\mathcal{Q}^{\nu,\ell}$ is a partition of $\mathbb{R}^2$.
However, at this point we cannot use the ``partition of
$\mathbb{R}^2$" on
$\sum_{\abs{\ell}\leq2^\nu}\sum_{P\in\mathcal{Q}^{\nu,\ell}}\tilde{\chi}_P$.
Therefore, if $0<q\leq 1$ we use the $q$-triangle inequality
$\abs{a+b}^q\leq \abs{a}^q+\abs{b}^q$ ($N>2/q$) or Hölder's
inequality if $1<q$ ($N>2$) to get (by hypothesis
$\lambda>3\max(1,1/q,1/p)>N$)
\begin{eqnarray*}
  \norm{f}_{\textbf{F}_1}
    &\lesssim& \norm{\left(\sum_{P\in\mathcal{Q}_{AB}}
        [2^{\nu(\alpha_1+1)}\abs{(s_{1,\lambda}^\ast)_P}\tilde{\chi}_P(\cdot)]^q\right)^{1/q}}_{L^p} \\
    &\lesssim& \norm{\left(\sum_{P\in\mathcal{Q}_{AB}}
        [\abs{P}^{-\alpha_2}\abs{(s_{1,\lambda}^\ast)_P}\tilde{\chi}_P(\cdot)]^q\right)^{1/q}}_{L^p}
        =\norm{\mathbf{s}_{1,\lambda}^\ast}_{\textbf{f}_2}.
\end{eqnarray*}
By Lemma \ref{l:s_ast_bnd_s} and Theorem \ref{t:S-T-psi-ops_Bnd} the
proof is complete.

\hfill $\blacksquare$ \vskip .5cm   

\subsection{Further relations}\label{sS:More_Relations_T-L-Shrl_T-L_Dyadic}
A dyadic cube at scale $\nu$ will be denoted by $Q_\nu$ and a shear
anisotropic ``cube" (a parallelepiped) at scale $j$ will be denoted
by $P_j$.
\begin{thm}\label{t:Dydc-T-L_fade_with_Shrlt-T-L}
Let $\alpha_1,\alpha_2\in \mathbb{R}$, $0<q_1,q_2\leq\infty$ and
$0<p_1,p_2<\infty$. Then, there exist sequences of functions
$\{f^{(j)}\}_{j\geq 0}$ such that
$\norm{f^{(j)}}_{\mathbf{F}^{\alpha_2,q_2}_{p_2}(AB)}\approx 1$, but
that $\norm{f^{(j)}}_{\mathbf{F}^{\alpha_1,q_1}_{p_1}}\rightarrow
0$, $j\rightarrow\infty$, when
$3(\alpha_2-1/p_2)>2\alpha_1-1/p_1+1$.
\end{thm}
\textbf{Proof}. For a sequence
$\mathbf{s}^{(j)}=\{s_{j,0,0}\}_{j\geq 0}$ such that
$\abs{s_{j,0,0}}=\abs{P_j}^{\alpha_2-\frac{1}{p_2}+\frac{1}{2}}$, we
have
$\norm{\mathbf{s}^{(j)}}_{\mathbf{f}^{\alpha_2,q_2}_{p_2}(AB)}=1$,
for all $j\geq 0$. Thus, $f^{(j)}(x)=s_{j,0,0}\psi_{j,0,0}(x)\in
\mathbf{F}^{\alpha_2,q_2}_{p_2}(AB)$ with
$\norm{f^{(j)}}_{\mathbf{F}^{\alpha_2,q_2}_{p_2}(AB)}\approx 1$.
From the compact support conditions on $\hat{\varphi}$ y
$\hat{\psi}$ the support of $(\psi_{j,0,0})^\wedge$ overlaps with
the support of $(\varphi_{\nu,0})^\wedge$ only when
$2j-5\leq\nu<2j$. Therefore, since assuming $\nu=2j$ and
$\abs{f^{(j)}\ast\varphi_{2^{2j}I}}=2^{2j}\abs{f^{(j)}\ast\varphi_{\nu,0}}$,
Lemma \ref{l:conv_shrlts-wvlts} yields
\begin{eqnarray*}
  \abs{f^{(j)}\ast\varphi_{2^{\nu}I}(x)}
    &=& \abs{\int_{\mathbb{R}^2} \abs{P_j}^{\alpha_2-\frac{1}{p_2}+\frac{1}{2}} \abs{\text{det }A}^{j/2}
        \psi(A^j(x-y))2^{4j}\varphi(2^{2j}y)dy} \\
    &\lesssim& 2^{-3j(\alpha_2-\frac{1}{p_2}+\frac{1}{2})+\frac{3j}{2}+4j}
        \int_{\mathbb{R}^2} \abs{\psi(A^j(x-y))}\abs{\varphi(2^{2j}y)}
        dy\\
    &\lesssim& \frac{2^{-3j(\alpha_2-\frac{1}{p_2})+j}}{(1+2^j\abs{x})^N},
\end{eqnarray*}
for every $N>2$. Then, for $j>0$ large enough and for $N$ such that
$Np_1>2$, we have
\begin{eqnarray*}
  \norm{f^{(j)}}_{\mathbf{F}^{\alpha_1,q_1}_{p_1}}
    &=& \norm{\left(\sum_{\nu=2j-5}^{2j-1} [2^{\nu\alpha_1}
        \abs{f^{(j)}\ast\varphi_{2^\nu I}}]^{q_1}\right)^{1/q_1}}_{L^{p_1}} \\
    &\leq& C_{N,q_1} \left(\int_{\mathbb{R}^2} 2^{2j\alpha_1p_1}\cdot
        \frac{[2^{-3j(\alpha_2-\frac{1}{p_2})+j}]^{p_1}}{(1+2^j\abs{x})^{Np_1}} dx\right)^{1/p_1} \\
    &=& C_{N,q_1}
        2^{2j\alpha_1-3j(\alpha_2-\frac{1}{p_2})+j-\frac{j}{p_1}},
\end{eqnarray*}
which tends to $0$ as $j\rightarrow\infty$ if
$2\alpha_1-\frac{1}{p_1}+1<3(\alpha_2-\frac{1}{p_2})$.

\hfill $\blacksquare$ \vskip .5cm   

\begin{thm}\label{t:ShrAnistr-T-L_in_ClssDydc-T-L}
Let $\alpha_1,\alpha_2\in \mathbb{R}$, $0<q_1,q_2\leq\infty$ and
$0<p_1,p_2<\infty$. Then, there exist sequences of functions
$\{f^{(\nu)}\}_{\nu\geq 0}$ such that
$\norm{f^{(\nu)}}_{\mathbf{F}^{\alpha_1,q_1}_{p_1}}\approx 1$, but
that $\norm{f^{(\nu)}}_{\mathbf{F}^{\alpha_2,q_2}_{p_2}(AB)}$
$\rightarrow 0$, $\nu\rightarrow\infty$, when
$2\alpha_1-4/p_1>3\alpha_2+1/q_2-1/p_2$.
\end{thm}
\textbf{Proof}. For a sequence
$\mathbf{s}^{(\nu)}=\{s_{\nu,0}\}_{j\geq 0}$ such that
$\abs{s_{\nu,0}}=\abs{Q_\nu}^{\frac{\alpha_1}{2}-\frac{1}{p_1}+\frac{1}{2}}$,
$\norm{\mathbf{s}^{(\nu)}}_{\mathbf{f}^{\alpha_1,q_1}_{p_1}}=1$, for
all $\nu\geq 0$. This means that
$f^{(\nu)}(x)=s_{\nu,0}\varphi_{\nu,0}(x)\in
\mathbf{F}^{\alpha_1,q_1}_{p_1}$ and
$\norm{f^{(\nu)}}_{\mathbf{F}^{\alpha_1,q_1}_{p_1}}\approx 1$, for
all $\nu\geq 0$. With the same arguments on the conditions of the
support of $(\varphi_{\nu,0})^\wedge$ and $(\psi_{j,\ell,0})^\wedge$
we assume $\nu=2j$ to get
\begin{eqnarray*}
  \abs{f^{(2j)}\ast\psi_{A^{-j}B^{-\ell}}(x)}
    &=& \abs{\int_{\mathbb{R}^2} s_{2j,0}2^{2j}\varphi(2^{2j}y)\abs{\text{det }A}^j\psi(B^\ell A^j(x-y)) dy} \\
    &\lesssim& 2^{-4j(\frac{\alpha_1}{2}-\frac{1}{p_1}+\frac{1}{2})+2j+3j}
        \int_{\mathbb{R}^2} \abs{\psi(B^\ell A^j(x-y))}\abs{\varphi(2^{2j}y)}dy \\
    &\lesssim&
       \frac{2^{-4j(\frac{\alpha_1}{2}-\frac{1}{p_1})}}{(1+2^j\abs{x})^N},
\end{eqnarray*}
for every $N>2$, by Lemma \ref{l:conv_shrlts-wvlts}. Then, for
$\nu>0$ large enough and $N>2$ such that $Np_2>2$,
\begin{eqnarray*}
  \norm{f^{(2j)}}_{\mathbf{F}_2}
    &\lesssim& \norm{\left(\sum_{\abs{\ell}\leq2^j} [\abs{P_j}^{-\alpha_2}\abs{f^{2j}\ast\psi_{A^{-j}B^{-\ell}}}]^{q_2}\right)^{1/q_2}}_{L^{p_2}} \\
    &\lesssim& \left(\int_{\mathbb{R}^2} (\sum_{\abs{\ell}\leq 2^j}
        [2^{3j\alpha_2}\frac{2^{-4j(\frac{\alpha_1}{2}-\frac{1}{p_1})}}{(1+2^j\abs{x})^N}]^{q_2})^{p_2/q_2} dx\right)^{1/p_2} \\
    &=& \left(\int_{\mathbb{R}^2} ((2^{j+1}+1)
        [\frac{2^{j(3\alpha_2-4(\frac{\alpha_1}{2}-\frac{1}{p_1}))}}{{(1+2^j\abs{x})^N}}]^{q_2})^{p_2/q_2} dx\right)^{1/p_2} \\
    &\leq& 2^{j(3\alpha_2-4(\frac{\alpha_1}{2}-\frac{1}{p_1})+\frac{1}{q_2})}
        \left(\int_{\mathbb{R}^2} \frac{1}{(1+2^j\abs{x})^{Np_2}} dx\right)^{1/p_2} \\
    &\lesssim& 2^{j(3\alpha_2-4(\frac{\alpha_1}{2}-\frac{1}{p_1})+\frac{1}{q_2})}\cdot
        2^{-\frac{j}{p_2}},
\end{eqnarray*}
which tends to $0$, as $j\rightarrow\infty$, if
$3\alpha_2-4(\frac{\alpha_1}{2}-\frac{1}{p_1})+\frac{1}{q_2}-\frac{1}{p_2}<0$.

\hfill $\blacksquare$ \vskip .5cm   

\vskip1cm
\section{Weights}\label{S:Wghts}
To extend this work to the weighted case $w\in
A_\infty=\cup_{p>1}A_p$, one can follow \cite{Bu82} and \cite{BH05}.
The spaces $\mathbf{F}^{\alpha,q}_p(AB)$ and
$\mathbf{f}^{\alpha,q}_p(AB)$ are then defined by $L^p(w)$
quasi-norms. For the weighted version of the Fefferman-Stein
inequality we refer the reader to \cite{AJ80} or to Remark 6.5 of
Chapter V in \cite{GC-RdF}. One adds $w\in A_{p_0}$ to the statement
of Lemma \ref{l:s_ast_bnd_s} and $N$ is chosen so that
$N>3\max(1,r/q,rp_0/p)$. Regarding the proof of Lemma
\ref{l:s_ast_bnd_s}, $\lambda$ should be chosen so that
$N>3\lambda/2>3\max(1,r/q,rp_0/p)$. For Theorem
\ref{t:S-T-psi-ops_Bnd} one adds $w\in A_\infty$. For its proof,
$\lambda$ should be chosen such that $0<\lambda<\min(p/p_0,q)$ and
$N$ such that $N>3\max(1,1/q,p_0/p)$.

\vskip1cm
\section{Proofs}\label{S:Proofs_lemmata}
We prove results of Subsections 3.1 and 4.2 in each of the
Subsections \ref{S:Proofs_lemmata}.$s$, $s=1,2$.

\vskip0.5cm
\subsection{Proofs for Subsection
\ref{sS:Notatn_&_Almost-Orthgnlty}}\label{sS:Proofs_Notation} In
order to prove Lemma \ref{l:c:l:conv_shrlts} we need two previous
results. The first one establishes that the sheared ellipsoids, in a
certain scale $j$ and for all shear parameter
$\ell=-2^j,\ldots,2^j$, contain the circle in a lower scale $j-1$.
This allows us to bound from below the distance of a point
$y\in\mathbb{R}^2$ under the $B^\ell A^j$ operation. The second one
is a basic ``almost orthogonality" result from which we derive Lemma
\ref{l:c:l:conv_shrlts}.

\begin{lem}\label{l:Nested_rectngls_ellipsoids}
Let $A$ and $B$ be as in Section \ref{S:Notation}. Then,
$$2^{j-1}\abs{x}<2^{j-\frac{1}{2}}\abs{x}< \abs{B^\ell A^j x},$$
for all $j\geq 0$, $\ell=-2^j,\ldots,2^j$ and all
$x\in\mathbb{R}^2$.
\end{lem}
\textbf{Proof}. For any $x\in\mathbb{R}^2$,
$$B^\ell A^j x =\left(\begin{array}{c}
              2^{2j}x_1+2^j\ell x_2 \\
              2^{j}x_2.
                               \end{array}
    \right)$$
With out loss of generality we can take $x\in\partial\mathbb{U}$,
where $\mathbb{U}$ is the unit disk. Then, $2^{j-1}\abs{x}=2^{j-1},$
and
\begin{eqnarray*}
  \abs{B^\ell A^jx}
    &=& ((2^{2j}x_1+2^j\ell x_2)^2+(2^jx_2)^2)^{1/2} \\
    &=& (2^{4j}x_1^2+2\cdot2^{3j}\ell x_1x_2+2^{2j}\ell^2x_2^2+(2^jx_2)^2)^{1/2} \\
    &=& (2^{4j}x_1^2+2\cdot2^{3j}\ell
        x_1(1-x_1^2)^{1/2}+2^{2j}\ell^2(1-x_1^2)+2^{2j}(1-x_1^2))^{1/2}.
\end{eqnarray*}
The extrema of $x_1(1-x_1^2)^{1/2}=\pm 1/2$ take place in
$x_1=\pm1/\sqrt{2}$ and $-2^{3j}\ell+2^{2j-1}\ell^2$ is minimized at
$\ell=2^j$. So,
$$\abs{B^\ell A^jx}\geq (2^{4j-1}-2^{3j}\ell+2^{2j-1}\ell^2+2^{2j-1})^{1/2}
    \geq 2^{j-\frac{1}{2}}.$$

\hfill $\blacksquare$ \vskip .5cm   

\noindent Similarly one can also prove that, for $i\geq j$,
$\abs{B^m A^iA^{-j}B^{-\ell} x}\geq 2^{i-j-1}\abs{x}$ holds true. To
see this observe that for $x\in\mathbb{U}$,
\begin{eqnarray*}
  \abs{B^mA^iA^{-j}B^{-\ell} x}
    &=& (2^{4(i-j)}x_1^2
        \pm2^{2(i-j)+1}(\ell 2^{2(i-j)}+m 2^{i-j})x_1(1-x_1^2)^{1/2} \\
    & &\;\;\;\; + (\ell 2^{2(i-j)}+m 2^{i-j})^2(1-x_1^2)
        + 2^{2(i-j)}(1-x_1^2))^{1/2}.
\end{eqnarray*}
With the same extrema $x_1(1-x_1^2)^{1/2}=\pm 1/2$ as before and
minimizing $-2^{2(i-j)}a+a^2/2$, where $a=\ell 2^{2(i-j)}+m
2^{i-j}$, the result follows as in the proof of Lemma
\ref{l:Nested_rectngls_ellipsoids}.

Observe that this result cannot be applied to the case of the
``discrete shearlets", since the shear parameter in this case runs
through $\mathbb{Z}$ making the ellipsoids thinner as
$\ell\rightarrow\pm\infty$ and consequently intersect the inner
circle.

The next result is fundamental to prove our first ``almost
orthogonality" result in form of a convolution in the space domain
(see (\ref{e:InnrProd_equiv_Conv_Shrlts})).
\begin{lem}\label{l:conv_shrlts}
Let $g,h\in\mathcal{S}$. Then, for all $N>2$ and $i=j-1, j, j+1\geq
0$, there exists $C_N>0$ such that
$$\abs{g_{j,\ell,k}\ast h_{i,m,n}(x)}\leq
\frac{C_N}{(1+2^j\abs{x-A^{-i}B^{-m}n-A^{-j}B^{-\ell}k})^{N}},$$ for
all $x\in\mathbb{R}^2$.
\end{lem}
\textbf{Proof.} First, since $g,h\in\mathcal{S}$ then $\abs{g(x)},
\abs{h(x)}\leq \frac{C_N}{(1+\abs{x})^N}$ for all $N\in\mathbb{N}$.
With the notation for $\varphi_M(x)$ we have
\begin{eqnarray*}
  g_{j,\ell,k}\ast h_{i,m,n}(x)
   &=& \int_{\mathbb{R}^2} 2^{3j/2}g(B^\ell A^jy-k)\cdot2^{3i/2}h(B^m A^i(x-y)-n) dy \\
   &=&  2^{-3(i-j)/2}g_{0,0,0}\ast h_{B^\ell A^jA^{-i}B^{-m}}(x'),
\end{eqnarray*}
with $x'=B^\ell A^jx-k-B^\ell A^jA^{-i}B^{-m}n=B^\ell
A^j(x-A^{-j}B^{-\ell}k-A^{-i}B^{-m}n)$. Following \cite[\S 6]{HW96},
define
\begin{eqnarray*}
  E_1 &=& \{y\in\mathbb{R}^2:\abs{x'-y}\leq 3\} \\
  E_2 &=& \{y\in\mathbb{R}^2:\abs{x'-y}> 3, \abs{y}\leq \abs{x'}/2\} \\
  E_3 &=& \{y\in\mathbb{R}^2:\abs{x'-y}> 3, \abs{y}> \abs{x'}/2\}.
\end{eqnarray*}
For $y\in E_1$ we have $1+\abs{x'}\leq 1+\abs{x'-y}+\abs{y}\leq
4(1+\abs{y})$. If $y\in E_3$, $1+\abs{x'}\leq 1+2\abs{y}\leq
2(1+\abs{y})$. When $y\in E_2$, $\frac{1}{2}\abs{x'}\leq
\abs{x'}-\abs{y}<\abs{x'-y}$, which implies
$4\abs{x'-y}=\abs{x'-y}+3\abs{x'-y}\geq 3+\frac{3}{2}\abs{x'}\geq 1
+\abs{x'}$ and  $1+\abs{B^m A^iA^{-j}B^{-\ell}(x'-y)}\geq
2^{i-j-1}\abs{x'-y}\geq 2^{i-j-3}(1+\abs{x'})$ (see Lemma
\ref{l:Nested_rectngls_ellipsoids}). Thus,
\begin{eqnarray*}
  \abs{g_{j,\ell,k}\ast h_{i,m,n}(x')}
    &\lesssim& \left\{\int_{E_1\cup E_3} + \int_{E_2} \right\}
        \frac{2^{-3(i-j)/2}}{(1+\abs{y})^{N}}
        \frac{\abs{\text{det }B^m A^iA^{-j}B^{-\ell}}}{(1+\abs{B^m A^iA^{-j}B^{-\ell}(x'-y)})^{N}} dy\\
    &\lesssim& \frac{2^{-3(i-j)/2}}{(1+\abs{x'})^{N}} \int_{E_1\cup E_3}
        \frac{\abs{\text{det }B^m A^iA^{-j}B^{-\ell}}}{(1+\abs{B^m A^iA^{-j}B^{-\ell}(x'-y)})^{N}}
        dy\\
    & & \;\;\;\; + \frac{2^{-3(i-j)/2}2^{3(i-j)}}{2^{(i-j-3)(N)}(1+\abs{x'})^{N}}
        \int_{E_2} \frac{1}{(1+\abs{y})^{N}}
        dy\\
    &\lesssim& \frac{C_N 2^{-3(i-j)/2}}{(1+\abs{x'})^{N}} \int_{\mathbb{R}^2}
        \frac{1}{(1+\abs{y})^{N}}
        dy \\
    & & \;\;\;\; + \frac{C_N2^{3(i-j)/2}}{2^{(i-j-3)(N)}(1+\abs{x'})^{N}}
        \int_{\mathbb{R}^2} \frac{1}{(1+\abs{y})^{N}}
        dy,
\end{eqnarray*}
for some $C_N>0$ and all $N>2$. The result follows from the fact
that $\abs{i-j}\leq 1$ and by replacing back $x'=B^\ell
A^j(x-A^{-j}B^{-\ell}k-A^{-i}B^{-m}n)$ in the estimates above and
because $\abs{B^\ell A^j y}> 2^{j-1}\abs{y}$ (see Lemma
\ref{l:Nested_rectngls_ellipsoids}) implies that
$$\frac{1}{(1+\abs{B^\ell A^j y})^N}<\frac{1}{(1+2^{j-1}\abs{y})^N}<
\frac{2^N}{(1+2^j\abs{y})^N}.$$

\hfill $\blacksquare$ \vskip .5cm   

As a corollary for Lemma \ref{l:conv_shrlts} we have our first
``almost orthogonality" property for the anisotropic and shear
operations for functions in $\mathcal{S}$.

\noindent\textbf{Proof of Lemma \ref{l:c:l:conv_shrlts}.} Identify
$(j,\ell,0)$ with $P$ and $(i,m,n)$ with $Q$. Write
$\abs{g_{A^{-j}B^{-\ell}}\ast
h_{i,m,n}(x)}=\abs{\abs{P}^{-1/2}g_P\ast h_Q}$. Since $\abs{i-j}\leq
1$, $\abs{P}^{-1/2}\thicksim\abs{Q}^{-1/2}$. Then, Lemma
\ref{l:conv_shrlts} yields
$$\abs{\abs{P}^{-1/2}g_P\ast h_Q}\leq \frac{C_N \abs{P}^{-1/2}}
    {(1+2^j\abs{x-x_Q})^N}\lesssim \frac{C_N \abs{Q}^{-1/2}}{(1+2^j\abs{x-x_Q})^N}.$$

\hfill $\blacksquare$ \vskip .5cm   

Our second ``almost orthogonality" result is stated in the Fourier
domain and gives more information since this time we take into
account the shear parameter $\ell$.

\noindent\textbf{Proof of Lemma \ref{l:Overlap_bnd}}. This is a
direct consequence of the construction and dilation of the
shearlets. Since $k$ and $n$ are translation parameters they do not
seem reflected in the support of $(\psi_{j,\ell,k})^\wedge$ or
$(\psi_{i,m,n})^\wedge$. By construction and by
(\ref{e:Discrt_Shrlt_Cond_Cone_1}) one scale $j$ intersects with
scales $j-1$ and $j+1$, only.

1) For one fixed scale $j$ and by (\ref{e:Discrt_Shrlt_Cond_Cone_2})
there exist \textbf{2} overlaps at the same scale $j$: one with
$(\psi_{j,\ell-1,k'})^\wedge$ and other with
$(\psi_{j,\ell+1,k''})^\wedge$ for all $k',k''\in\mathbb{Z}^2$.

2) Regarding scale $j-1$, one fixed $(\psi_{j,\ell,k})^\wedge$
overlaps with \textbf{3} other shearlets $(\psi_{j-1,m,k'})^\wedge$
at most for all $k,k'\in\mathbb{Z}^2$ because of 1) and because the
supports of the shearlets at scale $j-1$ have larger width than
those of scale $j$.

3) For a fixed scale $j$ consider the next three regions:
$\text{supp } (\psi_{j,\ell-1,k})^\wedge \cap \text{supp }
(\psi_{j,\ell,k'})^\wedge=R_{-1}$, $\text{supp }
(\psi_{j,\ell,k'})^\wedge \cap \text{supp }
(\psi_{j,\ell+1,k''})^\wedge=R_{+1}$ and $\text{supp
}(\psi_{j,\ell,k'})^\wedge \setminus (R_{-1}\cup R_{+1})=R_0$. Again
by construction, there can only be two overlaps for each $\xi$ at
any scale. Then, there exist at most two shearlets at scale $j+1$
that overlap with each of the three regions $R_i$, $i=-1,0,+1$ in
scale $j$: an aggregate of \textbf{6} for all translation parameters
$k,k',k''\in\mathbb{Z}^2$ at any scale $j$ or $j+1$.

Summing the number of overlaps at each scale gives the result.

\hfill $\blacksquare$ \vskip .5cm   

\vskip0.5cm
\subsection{Proofs of Subsection \ref{sS:Two_bsc_reslts}}\label{sS:Proof_Two_bsc_reslts}
To prove our results we follow \cite{FJ90}, \cite[\S 6.3]{HW96},
\cite{BH05} and \cite[\S 1.3]{Tr83}. Some previous well known
definitions and results are necessary.

\begin{defi}\label{d:Max_Fcn}
For a function $g$ defined on $\mathbb{R}^2$ and for a real number
$\lambda>0$ the \textbf{Peetre's maximal function} (see Lemma 2.1 in
\cite{Pe75}) is
$$g^\ast_\lambda(x)=\sup_{y\in\mathbb{R}^2}\frac{\abs{g(x-y)}}{(1+\abs{y})^{2\lambda}}, \;\;\; x\in\mathbb{R}^2.$$
\end{defi}

\begin{lem}\label{l:Bnd_diff_Max_Fcn}
Let $g\in \mathcal{S}'(\mathbb{R}^2)$ be such that $\text{supp
}(\hat{g})\subseteq \{\xi\in\hat{\mathbb{R}}^2: \abs{\xi}\leq R\}$
for some $R>0$. Then, for any real $\lambda>0$ there exists a
$C_\lambda>0$ such that, for $\abs{\alpha}=1$,
$$(\partial^\alpha g)^\ast_\lambda(x)\leq C_\lambda g^\ast_\lambda(x), \;\;\; x\in\mathbb{R}^2.$$
\end{lem}   
\textbf{Proof.} Since $g\in\mathcal{S}'$ has compact support in the
Fourier domain, $g$ is regular. More precisely, by the
Paley-Wiener-Schwartz theorem $g$ is slowly increasing (at most
polinomialy) and infinitely differentiable (e.g., Theorem 7.3.1 in
\cite{Ho90}). Let $\gamma$ be a function in the Schwartz class such
that $\hat{\gamma}(\xi)=1$ if $\abs{\xi}\leq R$. Then,
$\hat{\gamma}(\xi) \hat{g}(\xi) =\hat{g}(\xi)$ for all
$\xi\in\hat{\mathbb{R}}^2$. Hence, $\gamma\ast g=g$ and
$\partial^\alpha g=\partial^\alpha\gamma\ast g$. Moreover,
\begin{eqnarray*}
  \abs{\partial^\alpha g(x-y)}
    &=& \abs{\int_{\mathbb{R}^2} \partial^\alpha \gamma(x-y-z)g(z)dz}
        = \abs{\int_{\mathbb{R}^2} \partial^\alpha\gamma(w-y)g(x-w)dw} \\
    &\leq& \int_{\mathbb{R}^2} \abs{\partial^\alpha\gamma(w-y)}
    (1+\abs{w-y})^{2\lambda}
        (1+\abs{y})^{2\lambda} \frac{\abs{g(x-w)}}{(1+\abs{w})^{2\lambda}}dw,
\end{eqnarray*}
because of the triangular inequality. Therefore,
$$\abs{\partial^\alpha g(x-y)}\leq
g^\ast_\lambda(x)(1+\abs{y})^{2\lambda}
     \int_{\mathbb{R}^2} \abs{\partial^\alpha\gamma(w-y)} (1+\abs{w-y})^{2\lambda} dw.$$
Since $\gamma\in\mathcal{S}$, the last integral equals a finite
constant $c_\lambda$, independent of $y$, and we obtain
$$\abs{\partial^\alpha g(x-y)}\leq c_\lambda g^\ast_\lambda(x)(1+\abs{y})^{2\lambda},$$
which shows the desired result.

\hfill $\blacksquare$ \vskip .5cm   

We have a relation between the Hardy-Littlewood maximal function and
$g^\ast_\lambda$.
\begin{lem}\label{l:Max_Fcn_leq_HL-Max-Fcn}
Let $\lambda>0$ and $g\in\mathcal{S}'$ be such that $\text{supp
}(\hat{g})\subseteq \{\xi\in\hat{\mathbb{R}}^2: \abs{\xi}\leq R\}$
for some $R>0$. Then, there exists a constant $C_\lambda>0$ such
that
$$g^\ast_\lambda (x)\leq C_\lambda \left(\mathcal{M}(\abs{g}^{1/\lambda})(x)\right)^\lambda, \;\;\; x\in\mathbb{R}^2.$$
\end{lem}
\textbf{Proof.} Since $g$ is band-limited, $g$ is differentiable on
$\mathbb{R}^2$ (by the Paley-Wiener-Schwartz theorem), so we can
consider the pointwise values of $g$. Let $x,y\in\mathbb{R}^2$ and
$0<\delta<1$. Choose $z\in\mathbb{R}^2$ such that
$z\in\abs{B_\delta(x-y)}$. We apply the mean value theorem to $g$
and the endpoints $x-y$ and $z$ to get
$$\abs{g(x-y)}\leq\abs{g(z)}+\delta\sup_{w:w\in B_\delta(x-y)}
( \abs{\nabla g(w)}).$$ Taking the $(1/\lambda)^{\text{th}}$ power
and integrating with respect to the variable $z$ over
$B_\delta(x-y)$, we obtain
\begin{eqnarray}\label{e:l:Max_Fcn_leq_HL-Max-Fcn}
\nonumber
  \abs{g(x-y)}^{1/\lambda}
    &\leq& \frac{c_\lambda}{\abs{B_\delta(x-y)}}
      \int_{B_\delta(x-y)} \abs{g(z)}^{1/\lambda}dz \\
    &&+ c_\lambda \delta^{1/\lambda} \sup_{w:w\in B_\delta(x-y)}
    (\abs{\nabla g(w)})^{1/\lambda}.
\end{eqnarray}
Since $B_\delta(x-y)\subset B_{\delta+\abs{y}}(x)$,
$$\int_{B_\delta(x-y)} \abs{g(z)}^{1/\lambda}dz \leq \int_{B_{\delta+\abs{y}}(x)}
\abs{g(z)}^{1/\lambda}dz
    \leq \abs{B_{\delta+\abs{y}}(x)} \mathcal{M}(\abs{g}^{1/\lambda})(x),$$
and the $\sup$ term on the right hand side of
(\ref{e:l:Max_Fcn_leq_HL-Max-Fcn}) is bounded by
\begin{eqnarray*}
  \sup_{w:w\in B_{\delta+\abs{y}}(x)} (\abs{\nabla g(w)})^{1/\lambda}
    &=& \sup_{t:\abs{t}<\delta+\abs{y}} (\abs{\nabla g(x-t)})^{1/\lambda} \\
    &\leq& (1+\delta+\abs{y})^{2} \left[(\nabla g)^\ast_\lambda(x)\right]^{1/\lambda}.
\end{eqnarray*}
Substituting these last two inequalities in
(\ref{e:l:Max_Fcn_leq_HL-Max-Fcn}) yields
\begin{eqnarray*}
  \abs{g(x-y)}^{1/\lambda}
    &\leq& c_\lambda\frac{\abs{B_{\delta+\abs{y}}(x)}}{\abs{B_\delta(x-y)}}\mathcal{M}(\abs{g}^{1/\lambda})(x)  \\
    && + c_\lambda\delta^{1/\lambda}(1+\delta+\abs{y})^2 \left[(\nabla g)^\ast_\lambda(x)\right]^{1/\lambda},
\end{eqnarray*}
and since
$\abs{B_{\delta+\abs{y}}(x)}/\abs{B_\delta(x-y)}=(\delta+\abs{y})^2/\delta^2$,
we get
\begin{eqnarray*}
  \abs{g(x-y)}^{1/\lambda}
    &\leq& c_\lambda\frac{(\delta+\abs{y})^2}{\delta^2}\mathcal{M}(\abs{g}^{1/\lambda})(x)  \\
    && + c_\lambda\delta^{1/\lambda}(1+\delta+\abs{y})^2 \left[(\nabla g)^\ast_\lambda(x)\right]^{1/\lambda}.
\end{eqnarray*}
Taking the $\lambda^{\text{th}}$ power yields
$$\frac{\abs{g(x-y)}}{(1+\abs{y})^{2\lambda}}\leq
c'_\lambda\left\{\frac{1}{\delta^{2\lambda}}[\mathcal{M}(\abs{g}^{1/\lambda})(x)]^\lambda
    + \delta\left[(\nabla g)^\ast_\lambda(x)\right] \right\},$$
since $\delta<1$ implies $(1+\delta+\abs{y})\leq 2(1+\abs{y})$.
Taking $\delta$ small enough so that $c'_\lambda
C_\lambda\delta<1/4$ (where $C_\lambda$ is the constant in Lemma
\ref{l:Bnd_diff_Max_Fcn}) we obtain
$$g^\ast_\lambda(x)\leq c_\lambda\mathcal{M}(\abs{g}^{1/\lambda})(x)
    + \frac{1}{2}g^\ast_\lambda(x).$$
Assume for the moment that $g\in\mathcal{S}$, hence
$g^\ast_\lambda(x)<\infty$. So, we can subtract the second term in
the right-hand side of the previous inequality from the left-hand
side of the previous inequality and complete the proof for
$g\in\mathcal{S}$. To remove the assumption $g\in\mathcal{S}$, we
apply a standard regularization argument to a $g\in\mathcal{S}'$ as
done in p. 22 of \cite{Tr83} or in Lemma A.4 of \cite{FJ90}. Let
$\gamma\in\mathcal{S}$ satisfy $\text{supp }\hat{\gamma}\subset
B(0,1)$, $\hat{\gamma}(\xi)\geq 0$ and $\gamma(0)=1$. By Fourier
inversion $\abs{\gamma(x)}\leq 1$ for all $x\in\mathbb{R}^2$. For
$0<\delta<1$, let $g_\delta(x)=g(x)\gamma(\delta x)$. Then,
$\text{supp }\hat{g}_\delta$ is also compact,
$g_\delta\in\mathcal{S}$, $\abs{g_\delta}\leq \abs{g}$ for all
$x\in\mathbb{R}^2$ and $g_\delta\rightarrow g$ uniformly on compact
sets as $\delta\rightarrow 0$. Applying the previous result to
$g_\delta$ and letting $\delta\rightarrow 0$ we obtain the result
for general $\mathcal{S}'$.

\hfill $\blacksquare$ \vskip .5cm   

\noindent Lemmata \ref{l:Bnd_diff_Max_Fcn} and
\ref{l:Max_Fcn_leq_HL-Max-Fcn} are Peetre's inequality for
$f\in\mathcal{S}'$ whose proofs can be found in the references above
and we reproduce them for completeness.

\vskip0.5cm

\noindent\textbf{Proof of Lemma \ref{l:SupConv_leq_HLMaxFcn}.} Let
$g(x)=(\psi_{A^{-j}B^{-\ell}}\ast f)(x)$. Since $\psi$ is
band-limited, so is $g$. On one hand, since $j\geq 0$ and
$2^{j-1}\abs{y}\leq\abs{B^\ell A^j y}$,
\begin{eqnarray*}
  g^\ast_\lambda(t)
    &=& \sup_{y\in\mathbb{R}^2}\frac{\abs{g(t-y)}}{(1+\abs{y})^{2\lambda}}
        \geq \sup_{y\in\mathbb{R}^2} \frac{\abs{(\psi_{A^{-j}B^{-\ell}}\ast
            f)(t-y)}}{(1+2^j\abs{y})^{2\lambda}} \\
    &=& \sup_{y\in\mathbb{R}^2} \frac{\abs{(\psi_{A^{-j}B^{-\ell}}\ast
            f)(t-y)}}{2^{2\lambda}(2^{-1}+2^{j-1}\abs{y})^{2\lambda}} \\
    &\geq& 2^{-2\lambda} \sup_{y\in\mathbb{R}^2} \frac{\abs{(\psi_{A^{-j}B^{-\ell}}\ast
        f)(t-y)}}{(1+\abs{B^\ell A^jy})^{2\lambda}}
        = 2^{-2\lambda} \abs{(\psi_{j,\ell,\lambda}^{\ast\ast})(t)}.
\end{eqnarray*}
On the other hand,
\begin{eqnarray*}
  \mathcal{M}(\abs{g}^{1/\lambda})(t)
    &=& \sup_{r>0} \frac{1}{\abs{B_r(t)}} \int_{B_r(t)} \abs{(\psi_{A^{-j}B^{-\ell}}\ast f)(y)}^{1/\lambda} dy \\
    &=& \mathcal{M}(\abs{(\psi_{A^{-j}B^{-\ell}}\ast
        f}^{1/\lambda})(t).
\end{eqnarray*}
The result follows from Lemma \ref{l:Max_Fcn_leq_HL-Max-Fcn} with
$t=x$.
\hfill $\blacksquare$ \vskip0.5cm   


To prove Lemma \ref{l:s_ast_bnd_s} we need the next
\begin{lem}\label{l:bnd_conv_seqN-seq3-2}
Let $i\geq j\geq 0$ and $0<a\leq r$. Also, let $Q$ and $P$ be
identified with $(i,m,n)$ and $(j,\ell,k)$, respectively. Then, for
all $N>3r/a$, any sequence $\{s_P\}_{P\in\mathcal{Q}^{j,\ell}}$ of
complex numbers and any $x\in Q$,
$$(s_{r,N}^\ast)_Q:=\left(\sum_{P\in\mathcal{Q}^{j,\ell}}\frac{\abs{s_P}^r}{(1+2^j\abs{x_Q-x_P})^{N}}\right)^{1/r}
    \leq C_{a,r}
    \left[\mathcal{M}\left(\sum_{P\in\mathcal{Q}^{j,\ell}}\abs{s_P}^a\chi_P\right)(x)\right]^{1/a}.$$
Moreover, when $i=j$,
$$  \sum_{P\in\mathcal{Q}^{j,\ell}} \left[(s^\ast_{r,N})_P\tilde{\chi}_P(x)\right]^q
    \leq C_{a,r} \left[\mathcal{M}\left(\sum_{P\in\mathcal{Q}^{j,\ell}}
        (\abs{s_P}\tilde{\chi}_{P})^a\right)(x)\right]^{q/a}.$$
\end{lem}
\textbf{Proof.} Identify $(i,m,n)$ and $(j,\ell,k)$ with $Q$ and
$P$, respectively. Then, $x_Q=A^{-i}B^{-m}n$ and
$x_P=A^{-j}B^{-\ell}k$. Let
$\mathcal{Q}^{j,\ell}:=\{Q_{j,\ell,k}:k\in\mathbb{Z}^2\}$, then
$\mathcal{Q}^{j,\ell}$ is a partition of $\mathbb{R}^2$. Write
$d_P=\abs{x_Q-x_P}$. Thus, we bound the sum in the definition of
$(s_{r,N}^\ast)_Q$ as
$$\sum_{P\in\mathcal{Q}^{j,\ell}}\frac{\abs{s_P}^r}{(1+2^j\abs{x_Q-x_P})^N}
    \leq \left(\sum_{P\in\mathcal{Q}^{j,\ell}:d_P\leq 1} + \sum_{P\in\mathcal{Q}^{j,\ell}:d_P> 1}\right)
        \frac{\abs{s_P}^r}{(1+2^jd_P)^N}.$$
Choose $\lambda$ such that $N>3\lambda/2>3r/a$. Then, the inequality
$(2^j d_P)^N> (2^{3j/2}d_P)^\lambda$ holds whenever
$d_P>2^{j(3\lambda/2-N)/(N-\lambda)}$. So, the previous inequality
is bounded by
$$\sum_{P\in\mathcal{Q}^{j,\ell}:d_P\leq 1} \frac{\abs{s_P}^r}{(1+2^jd_P)^N}
    + \sum_{P\in\mathcal{Q}^{j,\ell}:d_P> 1} \frac{\abs{s_P}^r}{(2^{3j/2}d_P)^\lambda}.$$
Defining
\begin{eqnarray*}
  D_0
    &=& \{k\in\mathbb{Z}^2:\abs{A^{-i}B^{-m}n-A^{-j}B^{-\ell}k}\leq 1\} \\
    &=& \{P\in\mathcal{Q}^{j,\ell}:d_P=\abs{x_Q-x_P}\leq 1\}
\end{eqnarray*}
and
\begin{eqnarray*}
  D_\nu
    &=& \{k\in\mathbb{Z}^2:2^{\nu-1}<2^{3j/2}\abs{A^{-i}B^{-m}n-A^{-j}B^{-\ell}k}\leq 2^\nu\} \\
    &=& \{P\in\mathcal{Q}^{j,\ell}:2^{\nu-1}<2^{3j/2}d_P\leq 2^\nu\},\;\;\;
    \nu=1,2,3,\ldots,
\end{eqnarray*}
we have that
\begin{eqnarray*}
  \sum_{P\in\mathcal{Q}^{j,\ell}}\frac{\abs{s_P}^r}{(1+2^jd_P)^N}
    &\leq& \sum_{P\in D_0}\abs{s_P}^r +2^\lambda\sum_{\nu=1}^\infty\sum_{P\in D_\nu}\frac{\abs{s_P}^r}{2^{\nu\lambda}} \\
    &\leq& 2^\lambda\sum_{\nu=0}^\infty \sum_{P\in D_\nu} \frac{\abs{s_P}^r}{2^{\nu\lambda}}
        \leq 2^\lambda\sum_{\nu=0}^\infty 2^{-\nu\lambda}
            \left(\sum_{P\in D_\nu}\abs{s_P}^a\right)^{r/a},
\end{eqnarray*}
because $1+2^jd_P\geq 1$ and $a\leq r$. Now, when $x\in Q_{i,m,n}=Q$
and $P\in D_\nu$ then, by the definition of $D_\nu$,
$P=Q_{j,\ell,k}\subset B_{2^{\nu-3j/2+2}}(x)$
 (this holds because for $j\geq 0$ the
diameter of any $P=Q_{j,\ell,k}$ is less than $\sqrt{5}$ and the
intervals in the definition of the $D_\nu$'s are dyadic with
$\nu\geq 0$). Thus,
$$\sum_{P\in D_\nu}\abs{s_P}^a=\abs{P}^{-1}\int_{B_{2^{\nu-3j/2+2}}(x)}
    \sum_{P\in D_\nu}\abs{s_P}^a\chi_{P}(y)dy.$$
Hence, writing $\abs{\tilde{B}}=\abs{B_{2^{\nu-3j/2+1}}(x)}=\pi
2^{2\nu-3j+4}$ we have that for $x\in Q_{i,m,n}=Q$,
\begin{eqnarray*}
  \sum_{P\in\mathcal{Q}^{j,\ell}}\frac{\abs{s_P}^r}{(1+2^jd_P)^N}
    &\leq& 2^\lambda\sum_{\nu=0}^\infty 2^{-\nu\lambda}\left( \frac{\abs{P}^{-1}\abs{\tilde{B}}}{\abs{\tilde{B}}}
        \int_{\tilde{B}}\sum_{P\in D_\nu}\abs{s_P}^a\chi_{P}(y)dy \right)^{r/a} \\
    &\leq& 2^\lambda\sum_{\nu=0}^\infty 2^{-\nu\lambda}\pi^{r/a} 2^{4r/a}2^{2\nu r/a} \left(\mathcal{M}
        \left(\sum_{P\in D_\nu}\abs{s_P}^a\chi_{P}\right)(x)\right)^{r/a}\\
    &\leq& C_{r,a,\lambda}
        \left(\mathcal{M}\left(\sum_{P\in\mathcal{Q}^{j,\ell}}\abs{s_P}^a\chi_{P}\right)(x)\right)^{r/a},
\end{eqnarray*}
because $\abs{P}^{-1}=2^{3j}$ and $2r/a<\lambda$. 

To prove the second inequality multiply both sides by
$\tilde{\chi}_Q(x)$, rise to the power $q$ and sum over $Q\in
\mathcal{Q}^{j,\ell}$ to get
\begin{eqnarray*}
  \sum_{Q\in\mathcal{Q}^{j,\ell}} \left[(s^\ast_{r,N})_Q\tilde{\chi}_Q(x)\right]^q
    &\leq& C\sum_{Q\in\mathcal{Q}^{j,\ell}}  \left[\mathcal{M}\left(\sum_{P\in\mathcal{Q}^{j,\ell}}
        \abs{s_P}^a\chi_{P}\right)(x)\right]^{q/a}\tilde{\chi}_Q^q(x)\\
    &=& C\sum_{Q\in\mathcal{Q}^{j,\ell}}  \left[\mathcal{M}\left(\sum_{P\in\mathcal{Q}^{j,\ell}}
        (\abs{s_P}\tilde{\chi}_{P})^a\right)(x)\right]^{q/a}\chi_Q(x) \\
    &=& C \left[\mathcal{M}\left(\sum_{P\in\mathcal{Q}^{j,\ell}}
        (\abs{s_P}\tilde{\chi}_{P})^a\right)(x)\right]^{q/a},
\end{eqnarray*}
since $\mathcal{Q}^{j,\ell}$ is a partition of $\mathbb{R}^2$.

\hfill $\blacksquare$ \vskip .5cm   

\noindent\textbf{Proof of Lemma \ref{l:s_ast_bnd_s}}. Let $\lambda$
be such that $N>3\lambda/2>3\max(1,r/q,r/p)$. If $r<\min(q,p)$,
choose $a=r$. Otherwise, if $r\geq\min(q,p)$, choose $a$ such that
$r/(\lambda/2)<a<\min(r,q,p)$. It is always possible to choose such
an $a$ since $\lambda/2>\max(1,r/q,r/p)$ implies
$r/(\lambda/2)<\min(r,q,p)$. In both cases we have that
$$0<a\leq r<\infty, \;\;\; \lambda>2r/a, \;\;\; q/a>1, \;\;\; p/a>1.$$
The previous argument is similar to that of \cite{BH05}. Then, by
Lemma \ref{l:bnd_conv_seqN-seq3-2} and Theorem \ref{t:Feff-Stn_Ineq}
\begin{eqnarray*}
  \norm{\mathbf{s}^\ast_{r,N}}_{\mathbf{f}^{\alpha,q}_p(AB)}
    &=& \norm{\left(\sum_{P\in\mathcal{Q}_{AB}}(\abs{P}^{-\alpha}(s^\ast_{r,N})_P\tilde{\chi}_P)^q\right)^{1/q}}_{L^p} \\
    &=& \norm{\left(\sum_{j\geq 0}\sum_{\ell=-2^j}^{2^j}\abs{Q_{j,\ell}}^{-\alpha q}\sum_{P\in\mathcal{Q}^{j,\ell}}
        (s^\ast_{r,N})_P^q\tilde{\chi}_P^q \right)^{1/q}}_{L^p} \\
    &\leq& C\norm{\left(\sum_{j\geq 0}\sum_{\ell=-2^j}^{2^j}\abs{Q_{j,\ell}}^{-\alpha q}
        \left[\mathcal{M}\left(\sum_{P\in \mathcal{Q}^{j,\ell}}(\abs{s_P}\tilde{\chi}_{P})^a\right)\right]^{q/a} \right)^{1/q}}_{L^p} \\
    &=& C\norm{\left(\sum_{j\geq 0}\sum_{\ell=-2^j}^{2^j}
        \left[\mathcal{M}\left(\sum_{P\in \mathcal{Q}^{j,\ell}}(\abs{P}^{-\alpha}\abs{s_P}\tilde{\chi}_{P})^a\right)\right]^{q/a} \right)^{1/q}}_{L^p} \\
    &=& C\norm{\left(\sum_{j\geq 0}\sum_{\ell=-2^j}^{2^j}
        \left[\mathcal{M}\left(\sum_{P\in \mathcal{Q}^{j,\ell}}(\abs{P}^{-\alpha}\abs{s_P}\tilde{\chi}_{P})^a\right)\right]^{q/a}
        \right)^{a/q}}_{L^{p/a}}^{1/a}\\
    &\leq& C\norm{\left(\sum_{j\geq 0}\sum_{\ell=-2^j}^{2^j} \left(\sum_{P\in \mathcal{Q}^{j,\ell}}
        (\abs{P}^{-\alpha}\abs{s_P}\tilde{\chi}_P )^a\right)^{q/a}\right)^{a/q}}_{L^{p/a}}^{1/a} \\
    &=& C\norm{\left(\sum_{j\geq 0}\sum_{\ell=-2^j}^{2^j} \sum_{P\in \mathcal{Q}^{j,\ell}}
        (\abs{P}^{-\alpha}\abs{s_P}\tilde{\chi}_P )^q\right)^{a/q}}_{L^{p/a}}^{1/a} \\
    &=& C\norm{\left(\sum_{j\geq 0}\sum_{\ell=-2^j}^{2^j} \sum_{P\in \mathcal{Q}^{j,\ell}}
        (\abs{P}^{-\alpha}\abs{s_P}\tilde{\chi}_P )^q\right)^{1/q}}_{L^{p}} \\
    &=& C\norm{\mathbf{s}}_{f^{\alpha,q}_p(AB)},
\end{eqnarray*}
because $\mathcal{Q}^{j,\ell}$ is a partition of $\mathbb{R}^2$.

The reverse inequality is trivial since $\abs{s_Q}\leq
(s^\ast_{r,\lambda})_Q$ always holds.

\hfill $\blacksquare$ \vskip .5cm   

{\bf Acknowledgements}. The author thanks Keith Rogers for his
comments and questions on a very early draft. The author also thanks
Eugenio Hernández for his support and patience.


\end{document}